\documentclass{amsart}

\usepackage{amssymb}

\usepackage[all]{xy}

\newtheorem{thm}{Theorem}%[section]

\newtheorem{lem}[thm]{Lemma}
\newtheorem{cor}[thm]{Corollary}

\newtheorem{prop}[thm]{Proposition}

\newtheorem{conj}[thm]{Conjecture}
\theoremstyle{definition}
\newtheorem{defn}[thm]{Definition}

\newtheorem{say}[thm]{}
\newtheorem{exmp}[thm]{Example}

\newtheorem{prob}[thm]{Problem}
\newtheorem{ques}[thm]{Question}    %!!!!!!!!!!!!!!!!!!!!

\newtheorem{rem}[thm]{Remark}          

\newtheorem*{ack}{Acknowledgments}      % \renewcommand{\theack}{} 
\newtheorem{notation}[thm]{Notation}   
  
\newtheorem{defn-thm}[thm]{Definition--Theorem}  %!!!!!!!!!!!!!!!!!!!!!!!!
\newtheorem{defn-lem}[thm]{Definition--Lemma}  %!!!!!!!!!!!!!!!!!!!!!!!!
\theoremstyle{remark}

\let \cedilla =\c
\renewcommand{\c}[0]{{\mathbb C}}  

\renewcommand{\o}[0]{{\mathcal O}} 
\newcommand{\z}[0]{{\mathbb Z}}
\newcommand{\n}[0]{{\mathbb N}}
\renewcommand{\r}[0]{{\mathbb R}}

\newcommand{\s}[0]{{\mathbb S}}
\newcommand{\h}[0]{{\mathbb H}}
\newcommand{\p}[0]{{\mathbb P}}

\newcommand{\q}[0]{{\mathbb Q}}

\newcommand{\qtq}[1]{\quad\mbox{#1}\quad}
\newcommand{\spec}[0]{\operatorname{Spec}}

\newcommand{\mult}[0]{\operatorname{mult}}

\newcommand{\supp}[0]{\operatorname{Supp}}    
    
\newcommand{\codim}[0]{\operatorname{codim}}

\newcommand{\aut}[0]{\operatorname{Aut}}

\newcommand{\sing}[0]{\operatorname{Sing}}

\newcommand{\ord}[0]{\operatorname{ord}}

\newcommand{\tsum}[0]{\textstyle{\sum}} 
\newcommand{\tprod}[0]{\textstyle{\prod}}

\def\into{\DOTSB\lhook\joinrel\to}

\newcommand{\sV}{\mathcal{V}}
\newcommand{\sE}{\mathcal{E}}
\newcommand{\sW}{\mathcal{W}}
\newcommand{\sA}{\mathcal{A}}
\newcommand{\sC}{\mathcal{C}}

\newcommand{\g}{\gamma}
\def\loccoh#1.#2.#3.#4.{H^{#1}_{#2}(#3,#4)}

\DeclareMathAlphabet{\mathchanc}{OT1}{pzc}%
                                {m}{it}

\newcommand{\GL}{\mathrm{GL}}

\newcommand{\SL}{\mathrm{SL}}

\newcommand{\dd}[0]{{\mathbb D}}
\newcommand{\bdd}[0]{\overline{\mathbb D}}

\newcommand{\arc}[0]{\operatorname{Arc}}
\newcommand{\sharc}[0]{\operatorname{ShArc}}
\newcommand{\arco}[0]{\operatorname{Arc}^{\circ}}
\newcommand{\farc}[0]{\widehat{\operatorname{Arc}}}

\newcommand{\sharcr}[0]{\operatorname{ShArc}_{\r}}
\newcommand{\link}[0]{\operatorname{link}}
\newcommand{\hol}[0]{{\mathcal H}}

\begin{document}

\bibliographystyle{amsalpha}

%\today
\title{Holomorphic arcs on analytic  spaces}
\author{J\'anos Koll\'ar and Andr\'as N\'emethi}

\maketitle

The study of   formal arcs 
 was initiated by Nash in a 1967  preprint; published 
much later as \cite{nash-arc}. 
Arc spaces of smooth varieties have a rather transparent structure
but difficult problems arise for arcs passing through singularities.
The Nash conjecture on the irreducible components of such arc
spaces was proved for surfaces \cite{boba-pp} and for toric singularities 
\cite{MR2030097},
but counter examples were found in higher dimensions 
\cite{MR2030097, df-arc, k-nash2}.

Here we start the study of holomorphic arcs; these are
holomorphic maps $\phi:\bdd \to X$ of the closed unit disk to  a
complex analytic space $X$.
As one expects, there is not much conceptual difference between the set
of formal arcs and the set of holomorphic arcs since every formal arc
can be approximated by holomorphic arcs. However, a formal deformation of an
arc is a much more local object than a
  holomorphic deformation.
(See Remark \ref{rem.3.rem} and Example \ref{cusp.exmp.i} for more details.)
Thus, in many cases, the space of holomorphic arcs has infinitely many
connected components while the space of formal arcs always has only finitely
many.

For a  complex analytic space $X$ we define two variants --
the space of {\it arcs}, denoted by $\arc(X)$ and the
space of {\it short arcs}, denoted by $\sharc(X)$ --
and we study their connected components.
For short arcs we obtain complete answers for surface singularities
and  for  isolated quotient singularities in all dimensions.

\section{Arcs  on analytic  spaces}

\begin{defn}[Arcs]\label{arcs.defn}
Let $X$ be a  complex analytic space. 
Let $\dd\subset \c$ denote the open unit disk and $\bdd\subset \c$ its closure.
A  {\it  complex analytic arc} in $X$ is a holomorphic  morphism
$\phi:\bdd \to X$.
(That is, $\phi$ is defined and holomorphic  in some neighborhood of $\bdd$.)
The {\it center} of $\phi$ is the point $\phi(0)\in X$.

We think of an arc as describing the local behavior of a
morphism of  a Riemann surface to $X$.
The interesting aspects happen if the  image passes
through a singular point of $X$. We thus  localize at such a point
and get a complex analytic  morphism $\phi:\bdd \to X$
 such that $\supp\phi^{-1}(\sing X)=\{0\}$. This will be called
a {\it short  complex analytic arc} or a
{\it short arc.} 

More generally, if $Z\subset X$ is a subset then a {\it short arc}
on  $(Z\subset X)$ is a complex analytic  morphism $\phi:\bdd \to X$
 such that $\supp\phi^{-1}(Z)=\{0\}$. Thus a short arc
on $X$ is the same as a short arc on $(\sing X\subset X)$.

For many purposes it is easier to work with {\it  formal arcs.}
These are   morphisms
$\phi:  \spec \c[[t]]\to  X$.
Most earlier works on arc spaces considered formal arcs; see
\cite{nash-arc, MR1896234, MR2030097, MR2071659, MR2483946, boba-pp}
and the references there.
\end{defn}

\begin{defn}[Spaces of arcs]\label{fam.arcs.defn.1}
Let $X$ be a complex space and  $Z\subset X$  a subset. 
We consider various spaces of arcs.
First we define these only as sets, then we endow them with a natural topology.
\begin{enumerate}
\item $\arc(X)$ is the set of all arcs $\phi:\bdd\to X$.
\item $\arco(Z\subset X)$ is the set of those arcs $\phi$
for which  $\phi(\partial \bdd)\subset X\setminus Z$.
\item $\sharc(Z\subset X)$ is the set of short arcs on $Z\subset X$.
\end{enumerate}
\noindent If $Z=\sing X$ then we also use the  notation 
\begin{enumerate}\setcounter{enumi}{3}
\item $\arco(X):=\arco(\sing X\subset X)$ and
\item $\sharc(X):=\sharc(\sing X\subset X)$.
\end{enumerate}

Fixing a continuous metric  $d(\ , \ )$ on $X$, we get a
metric on $\arc(X)$ by setting
$$
d_a(\phi, \psi):= \sup\{d(\phi(t), \psi(t)): t\in \bdd\}.
$$
It is clear that the topology induced by this metric does not depend
on the choice of the  metric $d(\ , \ )$.
The other arc spaces  inherit their
topology from $\arc(X)$.

It is quite likely that these spaces also have  a natural structure as an
infinite dimensional complex space, but we do not establish this.
(See Conjecture \ref{local.atlas.conj} for a more precise version.) 
In all our proofs, we essentially write down finite dimensional
complex subspaces of $\arc(X)$ and work with them.

\end{defn}

\begin{rem}[Comparison of the arc spaces]\label{rem.3.rem}
The arc spaces $ \arc(X), \sharc(X)$ and the space of formal arcs,
which we denote by $\farc(X)$, are different as sets, but the main 
distinguishing feature comes from the deformations that we allow.

Consider for instance a short arc $\phi_0:(0\in \bdd)\to (Z\subset X)$
and a deformation of it as an arc  $\{\phi_s:s\in [0,1]\}$.
In general, $\phi_s^{-1}(Z)$ breaks up into several points;
write these as $p_1(s),\dots, p_m(s)\in \bdd$.

If $m>1$ then, from our point  of view,
$\phi_s$ is a global object. Picking any of the $p_i(s)$ and ignoring the
others corresponds to a family of formal arcs. Note that for each
$s\in [0,1]$ we can switch to a smaller disk $D_s:=\bdd(\epsilon_s)\ni p_i(s)$
of radius $\epsilon_s $
such that $\phi_s|_{D_s}$ becomes a short arc.
However, $\lim_{s\to 0}\epsilon_s=0$, thus we never get a
family of short arcs if $m>1$. 

If $m=1$ then, after a translation and rescaling,
we do get a family of short arcs.
Thus, working with short arcs is essentially an equisingularity
condition on families.

Spaces of formal arcs always have only finitely many
(connected or irreducible) components. Since an arc can be through of as
a collection of many formal arcs, it is not surprising that
$\arco(X)$ usually has infinitely many connected components.
By contrast, $\sharc(X)$ is essentially a dense subset of
$\farc(X)$. However, the equisingularity condition 
turns out to be quite restrictive, and frequently we again get
infinitely many connected components.
\end{rem}

\begin{exmp}\label{cusp.exmp.i} Consider the surface singularity
$$
S:=(xyz=x^4+y^4+z^4)\subset \c^3.
\eqno{(\ref{cusp.exmp.i}.1)}
$$
(This is one of the simplest   cusp singularities;
we study the general case in Sections \ref{sec.cusps}--\ref{sec.inoue}.)
Any arc  $\bigl(x(t), y(t), z(t)\bigr)$ on $S$ can be written
uniquely as
$$
\bigl(x(t), y(t), z(t)\bigr)=\bigl(t^mu(t), t^mv(t), t^mw(t)\bigr)
$$
where  $m\in \n$ and $\bigl(u(0), v(0), w(0)\bigr)\neq (0,0,0)$.
If $u(t)^4+v(t)^4+w(t)^4$ is nowhere zero on $\bdd$ then
 we can write this arc as
$$
\bigl(x(t), y(t), z(t)\bigr)=
\Bigl( u\cdot\frac{uvw}{u^4+v^4+w^4}, v\cdot\frac{uvw}{u^4+v^4+w^4}  ,
 w\cdot\frac{uvw}{u^4+v^4+w^4}  \Bigr).
\eqno{(\ref{cusp.exmp.i}.2)}
$$
Conversely (and this is a very special occurence)
for any $\bigl(u(t), v(t), w(t)\bigr)$ such that
 $$
u(t)^4+v(t)^4+w(t)^4\neq 0 \qtq{for} t\in \bdd,
\eqno{(\ref{cusp.exmp.i}.3)}
$$
  the formula (\ref{cusp.exmp.i}.2)
defines an arc on $S$.
(The condition (\ref{cusp.exmp.i}.3) is  open  on
$\arc(S)$ and it always holds on a smaller disc. It turns out that,
as far as  local deformations are concerned, the arcs satisfying 
(\ref{cusp.exmp.i}.3) are typical among all arcs.)
The center of this arc is the origin iff
$$
u(0)\cdot v(0)\cdot w(0)= 0.
\eqno{(\ref{cusp.exmp.i}.4)}
$$

Under the traditional, Nash version of the space of arcs
on $(0\in S)$, we can perturb $\bigl(u(t), v(t), w(t)\bigr)$
such that one of the functions vanishes at the origin with multiplicity
1 and the other two do not vanish. This shows that $\farc(0\in S)$
has 3 irreducible components, corresponding to which of the
three values $u(0),  v(0),  w(0)$ vanishes.

By contrast, $\sharc(0\in S)$ has infinitely many connected components.
Indeed, if one of the functions  $ u(t), v(t), w(t)$ vanishes
at a point  $t_0\in \dd$ then the arc (\ref{cusp.exmp.i}.2)
sends $t_0$ to the origin. Thus if the arc (\ref{cusp.exmp.i}.2) is  short
then the functions  $ u(t), v(t), w(t)$ vanish only at the origin.
This implies that the multiplicity of their zero at the origin is a 
locally constant function on $\sharc(0\in S)$.
(Strictly speaking, so far 
we have established this only for the open set where
(\ref{cusp.exmp.i}.3)  holds, but it turns out to be true in general.)
Thus we get 3 doubly infinite families of connected components
given by the possibilities
$$
\mult_0\bigl(u(t), v(t), w(t)\bigr)= (n_1, n_2, n_3)\in \n^3
\qtq{where} n_1 n_2 n_3=0.
$$

\end{exmp}

\begin{say}[Connected components of $\sharc$]\label{conn.comp.say}
Let $X$ be a complex space and $Z\subset X$ a closed subset such that
$X\setminus Z$ is connected. 
Given any  arc $\phi\in \arco(Z\subset X)$ we can restrict it to the boundary
to get
$$
\phi|_{\partial\bdd}: \partial\bdd \cong\s^1  \to X\setminus Z.
$$
There is no natural basepoint, thus this map defines an
element of 
$$
\pi_1\bigl(X\setminus Z\bigr)/(\mbox{\rm conjugation}).
$$
(We fix the counterclockwise orientation on $\partial\bdd $.)
The map is clearly locally constant  on $\arco(Z\subset X)$,
thus it descends to
$$
\pi_0\bigl(\arco(Z\subset X)\bigr)\to 
\pi_1\bigl(X\setminus Z\bigr)/(\mbox{\rm conjugation}).
\eqno{(\ref{conn.comp.say}.1)}
$$
Since $\sharc(Z\subset X)\subset \arco(Z\subset X)$,
we get similar maps for $\sharc(Z\subset X)$.

Note further that the semigroup  $(\r_{>0},+)$ acts freely on  
$\sharc(Z\subset X)$ by
$$
\sigma_s:\phi(t)\mapsto \phi\bigl(e^{-s}\cdot t\bigr)
$$
and we can view $\sharc(Z\subset X)$ as an interval bundle over the orbit 
space of this action. In particular, if $Z\subset U\subset X$ is any open 
subset then $\sharc(Z\subset X)$ and $\sharc(Z\subset U)$
are homeomorphic, so understanding the topology of
$\sharc(Z\subset X)$ is a local problem near $Z$.
(By contrast, $\arco(Z\subset U)$ does   depend on the choice of 
$Z\subset U\subset X$; see, however,  Conjecture \ref{indep.of.U.conj}.)

Thus, for an  isolated singularity $(0\in X)$,
the topology of $\sharc (0\in X)$ is independent of the
particular representative that we choose.
If $X$ is contractible then 
$X\setminus \{0\}$ is homotopy equivalent to
$\link(0\in X)$, thus we get a natural map
$$
w:\pi_0\bigl(\sharc(0\in X)\bigr)\to 
\pi_1\bigl(\link(0\in X)\bigr)/(\mbox{\rm conj}).
\eqno{(\ref{conn.comp.say}.2)}
$$
More generally, if $Z\subset X$ is compact,  
and $\link(Z\subset X) $ denotes the boundary of a 
regular neighborhood of $Z\subset X$ then
we have a natural map
$$
w:\pi_0\bigl(\sharc(Z\subset X)\bigr)\to 
\pi_1\bigl(\link(Z\subset X)\bigr)/(\mbox{\rm conj}).
\eqno{(\ref{conn.comp.say}.3)}
$$
We call  these the {\it  winding number maps.}
\end{say}

Our main result says that for surfaces these maps
 are injective.

\begin{thm}\label{sharc.surf.thm.2}
 Let $(0\in X)$ be a normal surface  singularity.
Then  the winding number map
(\ref{conn.comp.say}.2)
$$
\pi_0\bigl(\sharc(0\in X)\bigr)\to \pi_1\bigl(\link(0\in X)\bigr)/
(\mbox{\rm conj})
\eqno{(\ref{sharc.surf.thm.2}.1)}
$$
is an injection. 
\end{thm}

The map (\ref{sharc.surf.thm.2}.1) is rarely surjective since its image 
is clearly contained in 
$$
\ker\bigl[\pi_1\bigl(\link(0\in X)\bigr)\to \pi_1(E)\bigr]/
\bigl(\mbox{\rm conjugation by }\pi_1\bigl(\link(0\in X)\bigr)\bigr)
\eqno{(\ref{sharc.surf.thm.2}.2)}
$$
where $E$ is the exceptional curve of the minimal resolution of $(0\in X)$. 
The image of the map (\ref{sharc.surf.thm.2}.1)
can be explicitly given in terms of the minimal resolution, see
Theorem \ref{dlt.nonquot.thm}. 
The description shows that, with a few exceptions,
the image of the winding number map is  an infinite but  small
subset of  (\ref{sharc.surf.thm.2}.2).

Since the topology of the
minimal resolution is determined by the oriented homeomorphism type of
the link \cite{NP}, we obtain the following consequence.

\begin{cor} For a  normal surface  singularity  $(0\in X)$
the image of the winding number map (\ref{sharc.surf.thm.2}.1) --
and hence 
$\pi_0\bigl(\sharc(0\in X)\bigr)$ -- 
is determined by the oriented homeomorphism type of
$\link(0\in X) $. \qed
\end{cor}

For some interesting classes of surface singularites 
there are even better descriptions.  For
quotient singularities, this works in all dimensions.

\begin{thm}\label{sharc.quot.thm.1}
 Let $(0\in X)\cong (0\in \c^n)/G$ be an isolated  quotient  singularity.
Then there are natural identifications
$$
\pi_0\bigl(\sharc(0\in X)\bigr)= 
\pi_1\bigl(\link(0\in X)\bigr)/(\mbox{\rm conj})
= G/(\mbox{\rm conj}).
$$
\end{thm}

For $G\subset \SL(2,\c)$, this gives a concrete realization of the 
McKay correspondence between
nontrivial conjugacy classes of $G$ and
exceptional curves of the minimal resolution of $\c^2/G$.
We also get a McKay-type correspondence in higher dimensions, see
Paragraph \ref{mckay.say}. 

It turns out that $\pi_0\bigl(\sharc(0\in S)\bigr) $
is infinite for every other surface  singularity;  see
Theorem \ref{dlt.nonquot.thm}. 
We have especially clear  descriptions in two further cases.

Let   $(0\in S)$ be a surface cusp singularity.
(See  \cite{MR0393045} or Section \ref{sec.cusps} 
for the definition and basic properties.)
The fundamental group of the link of a cusp is given by an extension
$$
0\to H\cong \z^2 \to \pi_1\bigl(\link(0\in S)\bigr)\to \z \to 0.
$$
Under conjugation, $H$ acts trivially on itself, thus
conjugation by $\pi_1\bigl(\link(0\in S)\bigr) $
on $H$ can be described by a single matrix  $M\in \SL(2,\z)$
which can be easily computed  from the minimal resolution. 
The conjugacy classes are the orbits of the action generated by $M$.
The dual cusp  gives the same link but with a different orientation;
see Paragraph \ref{say.cusp.4}.

\begin{thm} \label{cusp.thm.maybe} Let $(0\in S)$ be a  surface 
cusp singularity.
Then
\begin{enumerate}
\item 
The image of the winding number  map
(\ref{conn.comp.say}.2) is contained in the subgroup $H$
modulo conjugation by $\pi_1\bigl(\link(0\in S)\bigr)$.
\item Every nontrivial  conjugacy class
of $H$ is obtained uniquely from a connected component of the space of
 short arcs in either the cusp
or the dual cusp for a suitable orientation of $\partial \bdd$. 
\end{enumerate}
\end{thm}

We get a different type of answer for 
normal surface singularities $(0\in S)$ with a good $\c^*$-action
given by $(\lambda, s)\mapsto \rho(\lambda,  s)$.
We say that an arc $\phi:\bdd\to (0\in S)$ is
{\it equivariant}  if there are $a, b\in \n$ such that 
$$
\phi\bigl(\lambda^b t\bigr)=\rho\bigl(\lambda^a, \phi(t)\bigr)
\qtq{for} t, \lambda\in \bdd.
$$ 
We can assume that $(a,b)=1$. Note that $b$ divides the order
of  the stabilizer subgroup of the $\c^*$-action along
the image of $\phi$, thus $b=1$ for most equivariant arcs.

\begin{prop} \label{quasihomog.prop} Let $(0\in X)$ be a 
normal surface singularity with a good $\c^*$-action.
\begin{enumerate}
\item Every connected component of $\sharc(0\in X)$ contains an
equivariant arc.
\item This equivariant arc is unique, except when 
\begin{enumerate}
\item either $(0\in X)$ is a quotient singularity
\item or $b=1$ in the notation above,
in which case there is an irreducible 1-parameter family.
\end{enumerate} 
\end{enumerate} 
\end{prop}

\begin{say}[The method of the proof of Theorem \ref{sharc.surf.thm.2}]
\label{method.sharc.surf.thm.2}
Let $(0\in X)$ be a normal surface singularity
and $f:Y\to X$ a proper bimeromorphic morphism
that is an isomorphism outside $0$. Set $Z:=\supp f^{-1}(0)$.
Composing with $f$ gives continuous bijections
$$
f_*:\arco(Z\subset Y)\to \arco(0\in X)
\qtq{and}
f_*:\sharc(Z\subset Y)\to \sharc(0\in X).
$$
However, these maps are usually {\bf not} homeomorphisms
and not even the number of connected components is preserved.

In general, there is a tension between two requirements.
\begin{itemize}
\item The local structure of
$(Z\subset Y) $ should be simple in order to be able to describe
$\sharc(Z\subset Y) $. 
\item The map $f$ should be simple so that
$f_*$ be a homeomorphism, or at least
that 
 $f_*:\pi_0\bigl(\sharc(Z\subset Y)\bigr)\to \pi_0\bigl(\sharc(0\in X)\bigr)
$
be a bijection.
\end{itemize}
Following the Nash conjecture \cite{nash-arc}, 
a first attempt could be to take $Y$ to be the minimal
resolution (or  the minimal log resolution) of $(0\in X)$.
Both of these fail the second requirement, even in very simple cases
like $X=(xy=z^n)$ for any  $n\geq 2$.

Our results  rest on the observation that
the {\it minimal dlt modification}
(to be defined in Section \ref{sec.dlt}) has both of the above good properties.
The proof uses  the tight connection between the
minimal dlt modification and the plumbing construction
of  links \cite{NP}. 

From the point of view of the plumbing construction, Theorems 
\ref{sharc.quot.thm.1}--\ref{cusp.thm.maybe} correspond to
the two exceptional cases and Proposition \ref{quasihomog.prop}
to the simplest case. 
\end{say}

\begin{say}[Description of the sections]
We start our discussions by describing the various arc spaces
of the unit disk $(0\in \dd)$ in Section \ref{sec.disk}. 
This is then used  in Section \ref{sec.snc} to determine
the connected components of the spaces of short arcs on
simple normal crossing pairs. We also study  equivariant versions.
Arcs on quotient singularities are related to equivariant arcs on
 simple normal crossing pairs.
This leads to the proof of Theorem \ref{sharc.quot.thm.1} in
Section \ref{sec.quot}. 

The Nash conjecture \cite{nash-arc}, proved in \cite{boba-pp} for surfaces,
describes the irreducible components of the space of formal
arcs on a surface $S$ in terms of the minimal resolution.
Similarly, in Section \ref{sec.dlt}, we use the 
{\it minimal dlt modification}
to understand the irreducible components of $\sharc(S)$.
This gives an
 explicit description of the image of the map (\ref{sharc.surf.thm.2}.1).

As with the Nash conjecture, it is easy to come up with a
list of candidates for the irreducible components of an arc space;
the hard part is to prove that none of them is contained in
the closure of another. In Section \ref{sec.fundgr}
we  study the conjugacy problem for the fundamental group
of links of surface singularities to prove that 
the map (\ref{sharc.surf.thm.2}.1) is  injective and
to give an explicit description of its image.
The particular cases of Seifert manifolds is discussed in 
(\ref{say:Seifert}--\ref{say:uj}).
At the end we also show Proposition \ref{quasihomog.prop}.

A topological proof of Theorem \ref{cusp.thm.maybe} is given in
Section \ref{sec.cusps}. A cusp and its dual can be seen very clearly
on hyperbolic Inoue surfaces. This leads to another proof 
in Section \ref{sec.inoue}.

Comments  and conjectures on long arcs,  higher dimensions and  real arcs are in
Section \ref{sec.conj}. 
\end{say}

 \begin{ack}
 We thank V.~Alexeev, 
T.~de~Fernex, F.~Forstneri{\v{c}}, H.~Hauser and P.~Popescu-Pampu
 for comments and references.
Partial financial support  to JK  was provided  by  the NSF under grant number 
DMS-07-58275 and by the Simons Foundation. 
Partial financial support  to AN  was provided  by OTKA Grants 81203 and 100796.
Part of this paper was written 
while JK  visited the University of Utah and Stanford University.
\end{ack}

\section{Families of arcs}

\begin{defn}[Families of arcs]\label{fam.arcs.defn}
Let $X$ be a complex space  and $V$ a topological space.
A {\it continuous family of arcs} parametrized by $V$ is a 
continuous function
$$
F(\ , \ ):V\times \bdd \to X\qtq{such that}
F(v, \ ): \bdd \to X \qtq{is holomorphic}
\eqno{(\ref{fam.arcs.defn}.1)}
$$
for every $v\in V$. 
If $V$ is a complex space and $F$ is holomorphic, we have a
{\it holomorphic family of arcs.}

A family of arcs is essentially the same as the corresponding
{\it classifying map}  $F_c:V\to \arc(X)$.

We have not defined a complex structure on $\arc(X)$,
but the following approximation allows us to think about
 holomorphic classifying maps.

Let $\hol(\bdd) $ denote the Banach space of holomorphic functions
of $\bdd$, that is, functions that are 
defined and holomorphic  in some neighborhood of $\bdd$.
An arc in $\c^N$ is the same as $N$ holomorphic functions
$$
\phi_1,\dots, \phi_N\in \hol(\bdd).
$$
Thus we can identify $\arc(\c^N)$ with
the complex Banach space $\hol(\bdd)^N$.
A norm is given by 
$$
d(\phi, \psi): \tsum_i \sup\{|\phi_i(t)- \psi_i(t)|: t\in \bdd\}.
$$
If $U\subset \c^N$ is an open set then $\arc(U)$ is naturally an open subset
of $\arc(\c^N)$ and if 
$X\subset U$ is a closed subset 
 defined by the equations
$f_i(x_1,\dots, x_N)=0$ then  $\arc(X)$ is a closed subset of $\arc(U)$
defined by the  equations
$$
f_i\bigl(\phi_1(t),\dots, \phi_N(t)\bigr)\equiv 0 \quad \forall\ i.
$$
Note that here we think of $f_i$ as a map
 $f_i:\arc(U)\to \hol(\bdd)$,
thus  $\arc(X)\subset \arc(U)$ has infinite codimension.

If $X$ is smooth of dimension $n$, then this construction
establishes a complex Banach manifold structure on $\arc(X)$,
locally modeled on $\hol(\bdd)^n$.

Note also that $\arc(\c^N)$ is a complex vector space, thus it
can be written as a union of finite dimensional complex vector spaces.
However, it seems that a typical $\arc(X)\subset \arc(\c^N)$
intersects all finite dimensional  vector subspaces
of $\arc(\c^N)$ in a finite set of points, thus we do not get
interesting analytic subspaces of  $\arc(X)$ this way.
\end{defn}

\begin{defn}[Irreducibility] \label{irred.defn} Let $X$ be a 
complex space and $A\subset \arc(X)$ a  subset.
 We say that $A$ is 
{\it strongly irreducible}
if for every finite subset $P\subset A$ there is a
(finite dimensional) irreducible complex space $V_P$ and a 
 holomorphic family 
of arcs $F_P:V_P\times \bdd\to X$ such that the image of the
classifying map $F_P:V_P\to \arc(X)$ lies in $A$ and contains $P$.

If $A$ is a (finite dimensional) complex space   then
$A$ is  strongly irreducible iff it is irreducible.

There are two simple cases when we prove that some
infinite dimensional arc space $A$ is strongly irreducible.

First, assume that we can write $A$ as a convex open subset
in some complex Banach space $\hol$. 
 If $V\subset \hol$ is any finite dimensional
vector subspace then $V\cap A$ is a convex open subset of
$V\cong \c^m$, hence an irreducible complex manifold.

Second,  assume that we can write $A$ as an open subset
of $M\times \hol$ where $M$ is a (finite dimensional) irreducible complex space
and $\hol$ is a complex Banach space. 
Assume further that $A$ contains a constant section $M\times \{v\}$ for some
$v\in \hol$ and  $\bigl(\{m\}\times \hol\bigr)\cap A$ is convex for every
$m\in M$. 
If $V\subset \hol$ is any finite dimensional
vector subspace containing $v$ then 
$A\cap  \bigl(M\times  V\bigr)$ is a finite dimensional complex manifold.
The fibers of its projection to $M $ are convex, open in $V\cong \c^m$
and $M\times \{v\}$ gives a section. Thus 
 $A\cap  \bigl(M\times  V\bigr)$ is irreducible.

Our definition is dictated by what we could prove in some cases and
it would be useful to develop a notion of irreducibility
of arc spaces that is more in line with the usual finite dimensional concept.
See Conjecture \ref{irred.prob} for further discussion of this topic.
\end{defn}

\section{Arc spaces of the disk}\label{sec.disk}

In this section 
we describe the connected components of the various arc spaces of the
disk $\dd$. We show that they are all complex Banach manifolds,
have finite dimensional approximations and irreducible in the strong sense
of (\ref{irred.defn}). In subsequent sections we build up
arc spaces of higher dimensional complex spaces from products
of  arc spaces of the disk.

\begin{say}[Arc spaces of the disk]\label{arc.disk.say.1}
One can think of 
$\arc(\dd)$ as the space of holomorphic functions
$\phi:\bdd\to \dd$. Thus $\arc(\dd)=B_{<1}\hol(\bdd)$, the open unit
ball in $\hol(\bdd)$, hence $\arc(\dd)$ is strongly irreducible
by (\ref{irred.defn}).

$\sharc(0\in \dd)$ is the space of holomorphic functions
$\phi:\bdd\to \dd$ whose only zero is at the origin.
The multiplicity of the zero is locally constant, giving
the connected components  $\{\sharc(0\in \dd)_m:m=1,2,\dots\}$.
It is convenient to set 
$\sharc(0\in \dd)_0:=\arc(\dd\setminus\{0\})$, the set of arcs
that do not pass through $0$. 

For $m=0$, the function $\phi$ has no zeros, thus it has a logarithm
$$
\log \phi: \bdd \to \{h\in\hol(\bdd): \Re (h)<0\}.
$$
The space of functions whose real part is everywhere negative is
a convex open subset of $\hol(\bdd)$, thus has the same
type of finite dimensional approximations as the unit ball.
The actual arc space is the quotient
$$
\sharc(0\in \dd)_0\cong \{h\in\hol(\bdd): \Re (h)<0\}/2\pi i \z.
$$

If $\phi$ has an $m$-fold zero at the origin
then $t^{-m}\phi$ has no zeros. On the boundary of $\bdd$ it is still
strictly less than 1 in absolute value, hence, by the maximum principle,
 $t^{-m}\phi\in \sharc(0\in \dd)_0$. Thus multiplication by $t^m$ gives an
isomorphism
$$
t^m: \sharc(0\in \dd)_0\cong \sharc(0\in \dd)_m.
$$

$\arco(0\in \dd)$ is the space of holomorphic functions
$\phi:\bdd\to \dd$ that have no zero on the boundary.
The number  of  zeros in $\dd$ is locally constant, giving
the connected components  $\{\arco(0\in \dd)_m:m\in \n\}$.
Note that $\arco(0\in \dd)_0=\sharc(0\in \dd)_0$.

For $m>0$ we can use  Blaschke products to write any
$\phi\in \arco(0\in \dd)_m$ uniquely in the form
$$
\phi(t)=g(t)\cdot \tprod_{i=1}^m \frac{t-a_i}{1-\bar a_i t}
$$
where $g\in \arco(0\in \dd)_0$ and $a_i$ are the zeros of $\phi$. 
This gives a surjective real analytic map
$$
\dd^m\times \{h\in\hol(\bdd): \Re (h)<0\}\to \arco(0\in \dd)_m
$$
which is equivariant with respect to the
$2\pi i \z$ action on the second factor and the
symmetric group action permuting the coordinates on the first factor.
This shows that $\arco(0\in \dd)_m$ is connected, but it does not
give enough complex analytic subspaces.

In order to get the complex structure right, we write
$$
\phi(t)=\psi(t)\cdot \textstyle{\tprod}_{i=1}^m (t-a_i)
$$
where $a_i$ are the zeros of $\phi$.
The condition on $\psi(t)$ is more complicated than before; we need that
$$
0<|\psi(t)|<\Bigl|\tprod_{i=1}^m \frac{1}{t-a_i}\Bigr|\qtq{for every} t\in \bdd.
$$
As before, by passing to $\log \psi$   this becomes a  convex condition
$$
\Re \bigl(\log \psi(t)\bigr)<
\log \Bigl|\tprod_{i=1}^m \frac{1}{t-a_i}\Bigr|\qtq{for every} t\in \bdd.
$$
 Thus we get a
holomorphic parametrization of $\arco(0\in \dd)_m$
by an open subset 
$$
U_m\subset \dd^m\times  \hol(\bdd)
$$
whose intersection with  each 
$\{(a_1,\dots, a_m)\}\times \hol(\bdd)$ is convex.
Note that $U_m$ does not contain any constant section  $\dd^m\times \{c\}$,
but if $\epsilon>0$ and $\Re(c)$ is very negative then 
$\dd(1-\epsilon)^m \times \{c\}\subset U_m$. 
This is enough to show that $U_m$ is strongly irreducible
by (\ref{irred.defn}).
 \end{say}

In Section \ref{sec.quot}, we will also need equivariant versions of the
above arc spaces.

\begin{say}[Equivariant arc spaces of the disk]\label{arc.disk.say}
Let $G\subset \aut(0\in \bdd)$ be a finite subgroup and
$\rho:G\to \aut(Z\subset X)$ a representation. An arc $\phi$ is called
{\it $\rho$-equivariant} if 
$$
\phi\bigl(g(t)\bigr)=\rho(g)\bigl(\phi( t)\bigr).
$$
The set of $\rho$-equivariant short arcs is denoted by
$\sharc(Z\subset X)^{\rho}$.

Note that $\aut(0\in \bdd)$ is just the group of rotations,
thus every finite subgroup is cyclic and generated by
$\epsilon=e^{2\pi i/m}$ for some $m\in \n$.

For arc spaces of a disk there are even fewer possibilities.

Fix  natural numbers $m>0, a\geq 0$ and a primitive $m$th root of unity
$\epsilon$.

Let $\bdd(t)$ be a closed disk with coordinate $t$ and $\z/m$-action
$t\mapsto \epsilon t$  and $\dd(z)$ be an open disk with coordinate $z$ and 
$\z/m$-action
$z\mapsto \epsilon^a z$. 
The representation $\rho=\rho(a)$ is determined by
 $a$ modulo $m$.

We study the spaces of   $\rho(a)$-equivariant arcs.
These can be thought of as  functions $z=\phi(t)$ such that
$$
\phi(\epsilon t)=\epsilon^a\phi(t).
\eqno{(\ref{arc.disk.say}.1)}
$$
Every such $\phi$ can be uniquely written as
$\phi(t)=t^a\psi(t^m)$ for some $\psi$ in the  unit ball of
$\hol(\bdd(s))$ where $s=t^m$. Thus the space of  $\rho(a)$-equivariant arcs
is connected and strongly irreducible.

The space  $\sharc(0\in \dd(z))_r$ contains a $\rho(a)$-equivariant arc
only if $r\equiv a\mod m$ and then these can be written as
$\phi(t)=t^r\psi(t^m)$ for some invertible $\psi$ in the  unit ball of
$\hol(\bdd(s))$ where $s=t^m$.  Thus 
multiplication by $z^r$ identifies
$\sharc_0(\dd)_0=\arc(\dd\setminus\{0\})$ 
(where we think of $\dd$ as the unit disk with coordinate $s=t^m$) with 
the $\rho(a)$-equivariant arcs in $\sharc(0\in \dd)_r$.
\end{say}

\section{Arc spaces of the polydisk}\label{sec.polydisk}

It is clear from the definition that
$$
\arc\bigl(\textstyle{\prod}_{i\in I} X_i\bigr)=
\textstyle{\prod}_{i\in I}\arc(X_i).
$$
In particular, $\arc(\dd^n)\cong \arc(\dd)^n$
is a convex open set in the Banach space
$\hol(\bdd)^n$. 

By contrast, the spaces of short arcs have a more complicated
behavior. To start with, there are several ways to define
the product of pairs $(Z_i\subset X_i)$.

\begin{defn}[Products of pairs]\label{prod.pair.defn}
For simple normal crossing pairs the natural notion is
 the  product (or maximal product) of pairs $(Z_i\subset X_i)$, defined as
$$
\tprod_{i\in I}(Z_i\subset X_i):=
\Bigl(\bigcup_{i\in I}
Z_i\times \textstyle{\prod}_{j\neq i}X_j\subset 
\textstyle{\prod}_{i\in I}X_i\Bigr).
\eqno{(\ref{prod.pair.defn}.1)}
$$
The minimal product is
$$
\Bigl(\textstyle{\prod}_{i\in I}Z_i\subset\textstyle{\prod}_{i\in I} X_i\Bigr).
\eqno{(\ref{prod.pair.defn}.2)}
$$
\end{defn}

\begin{say}[Short arcs on maximal products]\label{snc.short.arcs.say}
In order to write the space of short arcs on a product
it is convenient to set
$$
\sharc^+(Z\subset X):=\sharc(Z\subset X)\textstyle{\coprod} \arc(X\setminus Z).
\eqno{(\ref{snc.short.arcs.say}.1)}
$$
Then
$$
\sharc^+\Bigl(\textstyle{\prod}_{i\in I}(Z_i\subset X_i)\Bigr)=
\prod_{i\in I}\sharc^+\bigl(Z_i\subset X_i\bigr).
\eqno{(\ref{snc.short.arcs.say}.2)}
$$
In particular,
$$
\sharc^+\bigl((0\in \dd)^r\times \dd^s\bigr)=
\bigl(\sharc^+(0\in \dd)\bigr)^r\times \bigl(\arc(\dd)\bigr)^s.
\eqno{(\ref{snc.short.arcs.say}.3)}
$$
The latter has  an obvious equivariant version.
If $\z/m$ acts on $\bdd$ by $t\mapsto \epsilon t$ and on
the target by
$$
\rho:(x_1,\dots, x_r, y_1,\dots, y_s)\mapsto
 (\epsilon^{a_1}x_1,\dots, \epsilon^{a_r}x_r,
\epsilon^{b_1} y_1,\dots, \epsilon^{b_s}y_s)
$$
then
$$
\sharc^+\bigl((0\in \dd)^r\times \dd^s\bigr)^{\rho}=
\tprod_{i=1}^r\sharc^+(0\in \dd)^{\rho(a_i)}\times 
\tprod_{i=1}^s\arc(\dd)^{\rho(b_i)}.
\eqno{(\ref{snc.short.arcs.say}.4)}
$$
\end{say}

It is harder to study short arcs on a minimal product
since the product of two non-short arcs is frequently short.
Especially for singular spaces, the precise answer seems complicated.

In our applications the ambient space is  
the polydisc but the subvariety sits between the
minimal and the maximal   product.

\begin{defn} \label{Z2.defn} 
Let  $\arc(0\in \dd)\subset \arc(\dd)$ denote the set of those arcs $\phi$
for which $\phi(0)=0$. 
Let $Z_2\subset \dd^n$ be the union of the
codimension 2 coordinate hyperplanes $\bigcup_{i\neq j}(x_i=x_j=0)$. 
Set
$$
\sharc^*(0\in \dd^n):=\bigl\{\phi\in \arc(0\in \dd)^n: 
\supp\phi^{-1}(Z_2)=\{0\}\bigr\}.
$$
It is clear that
$$
\sharc^*(0\in \dd^n)\subset \sharc(0\in \dd^n)\subset\arc(0\in \dd)^n.
$$
More generally, given representations $\rho_i:\z/m\to \c^*$
and their product $\rho:\z/m\to (\c^*)^n$ set
$$
\sharc^*(0\in \dd^n)^{\rho}
:=\Bigl\{\phi\in \sharc(0\in \dd^n)^{\rho}: 
\supp\phi^{-1}(Z_2)=\{0\}\Bigr\}.
$$
\end{defn}

\begin{prop}\label{shrac*.prop}  
For $n\geq 2$ the arc space
$\sharc^*(0\in \dd^n)$ is connected, strongly irreducible
and dense in $\arc(0\in \dd)^n $. Thus 
$\sharc(0\in \dd^n)$ is also
connected.

The same holds for 
 $\sharc^*(0\in \dd^n)^{\rho}$ and $\sharc(0\in \dd^n)^{\rho}$.
\end{prop}

Proof. For $i\neq j$, let $W_{ij}\subset \arc(0\in \dd)^n$
be the set of those functions $(f_1,\dots, f_n)$ such that
$f_i$ amd $f_j$ have a common zero in $\dd\setminus \{0\}$.
These  exactly correspond to arc that intersect the
codimension 2 coordinate hyperplane $(x_i=x_j=0)$ outside the origin. Thus
$$
\sharc^*(0\in \dd^n)=\arc(0\in \dd)^n\setminus \bigcup_{i\neq j} W_{ij}.
$$
The structure of the subsets $W_{ij}$ is somewhat complicated
(they are only real subanalytic), but
we show in
(\ref{com.zeros.say})
  that they are locally contained in a complex hypersurface.

By Lemma \ref{remove.codim.2.lem}, this implies that 
 $\sharc^*(0\in \dd^n)$
is connected, strongly irreducible and dense in $\sharc(0\in \dd)^n $. 
Thus any subspace in between these two is also connected, hence
$\sharc(0\in \dd^n)$ is connected.

The argument is the same for the equivariant versions.\qed
\medskip

More generally,  we study the locus where a collection of
holomorphic functions has an unexpected common zero.

\begin{say}[Common zeros of holomorphic functions]\label{com.zeros.say}
 Let $V_1,\dots, V_k$ be irreducible complex spaces
parametrizing holomorphic functions 
$F_i:V_i\times \bdd\to \dd$. For $f_1,\dots , f_k\in \hol(\bdd)$ let
$Z\bigl(f_1,\dots , f_k\bigr)\subset \bdd$ denote their common zero set
(with multiplicity).
The   common zero set of the families $F_1,\dots, F_k$ is then 
$$
Z\bigl(F_1,\dots, F_k\bigr)=
\bigcap  \bigl\{Z\bigl(f_1,\dots , f_k\bigr) : 
(f_1,\dots, f_k)\in V_1\times\dots\times V_k\bigr\}
$$
where we again  keep track of the multiplicities.
Let 
$$
W^{\rm exc}=W^{\rm exc}\bigl(F_1,\dots, F_k\bigr)
\subset V_1\times\dots\times V_k
$$
denote the set of those $(f_1,\dots, f_k) $ 
that have extra zeros, that is
for which
$$
Z\bigl(f_1,\dots , f_k\bigr)\supsetneq Z\bigl(V_1,\dots, V_k\bigr).
$$
$W^{\rm exc}$ is a subset of $V_1\times\dots\times V_k$,
which can be understood as follows.

Pick any function $g\in \hol(\bdd)$ whose zero set equals
$Z\bigl(V_1,\dots, V_k\bigr)$. Replacing each $F_i$ by $F_i/g$,
we get new families of holomorphic functions on $\bdd$
such that
$$
Z\bigl(F_1/g,\dots, F_k/g\bigr)=\emptyset.
$$
Furthermore,
$$
W^{\rm exc}\bigl(F_1,\dots, F_k\bigr)=
W^{\rm exc}\bigl(F_1/g,\dots, F_k/g\bigr).
$$
Thus, it is sufficient to understand $W^{\rm exc}$
in the special case when $Z\bigl(F_1,\dots, F_k\bigr)=\emptyset$.
We assume the latter from now on.

Let $H_i\subset V_1\times\dots\times V_k\times \bdd$
be the pull-back of the zero set of $F_i$
by the $i$th coordinate projection. Then
$W^{\rm exc}$ is the image of $H_1\cap\cdots\cap H_k$
under the projection 
$$
\Pi: V_1\times\dots\times V_k\times \bdd\to V_1\times\dots\times V_k.
$$
$\Pi$ is proper since $\bdd$ is compact,  so
$W^{\rm exc}$ is closed in $V_1\times\dots\times V_k$.

To understand its local structure, pick a 
 point $ (f_1,\dots, f_k)\in V_1\times\dots\times V_k$
such that not all the $f_i$ are identically 0.
We can extend the functions $F_i$ to a larger open disk
$\dd(1+\epsilon)$ where some of the $f_i$ has only finitely many zeros.
Then $\Pi$ restricts to a finite morphism $\Pi(\epsilon)$ on
$$
H_1(\epsilon)\cap\cdots\cap H_k(\epsilon)
\subset V_1\times\dots\times V_k\times \bdd(1+\epsilon).
$$
The image of $\Pi(\epsilon)$ describes those  $ (f_1,\dots, f_k)$
that have a common zero in $\dd(1+\epsilon)$. Thus
$W^{\rm exc}$ is a subset of this. 

Of course this is  useful information only if
$\Pi(\epsilon)$ is not dominant near  $ (f_1,\dots, f_k)$,
that is when $H_1(\epsilon)\cap\cdots\cap H_k(\epsilon) $
has codimension $\geq 2$ near
$ \{(f_1,\dots, f_k)\}\times \bdd(1+\epsilon)$.

We claim that this holds if at least 2 of the  $f_i$ are not identically 0.

To see this, pick a point $z_0\in Z(f_1,\dots, f_k)$. Since
$Z\bigl(F_1,\dots, F_k\bigr)=\emptyset$, there is at least one index,
say $i=1$, and a 1-parameter family  $f_1(s)\subset V_1$ such that
$f_1(s)(z_0)\neq 0$ for $s\neq 0$. Since 
at least 2 of the  $f_i$ are not identically 0, up to re-indexing we may 
assume that $f_2$ is not identically 0. 
Near $z_0$ the only zero of $f_2$ is $z_0$ which is
not a zero of $f_1(s)$ for $s\neq 0$.
Thus the restriction of $\Pi$ to 
$H_1\cap H_2$ is not dominant near 
$ \{(f_1,\dots, f_k)\}\times\{z_0\}$, hence
$H_1\cap H_2$ 
has codimension $\geq 2$ near
$ \{(f_1,\dots, f_k)\}\times\{z_0\}$. 
\end{say}

\begin{defn} Let $V$ be a complex space
and $W\subset V$ a closed subset. We say that
$W$ has {\it complex codimension $\geq k$}
if each point $w\in W$ has an open neighborhood $w\in U_w\subset V$
and a closed, complex subspace  $Z_w\subset U_w$ of codimension $\geq k$
such that $W\cap U_w\subset Z_w$.
We are mostly interested in the case $k=1$.

A finite union of  closed subsets of complex codimension $\geq 1$
also has complex codimension $\geq 1$.
\end{defn}

\begin{lem} \label{remove.codim.2.lem}
Let $V$ be a normal,  irreducible complex space
and $W\subset V$ a closed subset of complex codimension $\geq 1$.
Then $V\setminus W$ is also an  irreducible complex space.
\end{lem}

Proof. Let $Y\subset V\setminus W$ be a  irreducible component.
Let $w\in W$ be any point. Then $Y\cap (U_w\setminus Z_w)$
is an irreducible component of  $U_w\setminus Z_w $.
Since $U_w$ is normal, the latter is irreducible.
Thus if $Y\cap (U_w\setminus Z_w)\neq\emptyset$ then
the closure $\bar Y\subset V$ contain $U_w$. Thus
$\bar Y$ is a closed analytic subset of $V$, hence an
irreducible component. Thus $\bar Y=V$ and so
$Y=V\setminus W$.\qed

\medskip

By applying the above arguments to suitable projections
we obtain the following result, which implies that
we can move arcs away from codimension 2 subsets of the
smooth locus without changing the connected components
of the arc spaces.

\begin{prop}\label{codim2.remove.prop} Let $X$ be a complex space and
$Z\subset Z_2\subset X$  closed subsets. 
Assume that $X$ is smooth along $Z_2\setminus Z$ and 
$\codim_X(Z_2\setminus Z)\geq 2$.  Then
$$
\bigl\{\phi\in \sharc(Z\subset X):
\phi^{-1}(Z_2\setminus Z)\neq \emptyset\bigr\}
$$ has  
complex codimension $\geq 1$ in $\sharc(Z\subset X) $. \qed
\end{prop}

\section{Simple normal crossing pairs}\label{sec.snc}

\begin{say}[Intersection number of arcs and divisors]\label{snc.short.arcs.say.2}
Let $Y$ be a complex manifold and 
$D\subset Y$ a divisor. For some $Z\subset D$, let
$\phi\in \arco(Z\subset X)$ be an arc. We claim that
the intersection number $(\phi\cdot D)$ is defined and that it is a
locally constant function on  $\arco(Z\subset X)$.

The intersection number will be a sum of terms
for each point  $p\in \supp \phi^{-1}(D)$. To define the local contribution
at $p$, we restrict $\phi$ to a smaller disk $\dd(p,\epsilon)$ around $p$.
Let $\phi(p)\in U_p\subset X$ be an open neighborhood
such that $D\cap U_p=(F_p=0)$ for some holomorphic function $F_p$ on $U_p$.
By choosing  $\epsilon$ small enough, we may assume that
$\phi\bigl(\dd(p,\epsilon)\bigr)\subset U_p$. Then
$F_p\circ \phi$ is a holomorphic function on $\dd(p,\epsilon)$
and the local contribution is the multiplicity of its zero at $p$.

If we perturb $\phi$, these numbers stay constant as long as
the image of $\partial \dd(p,\epsilon)$ stays disjoint from $D$.
Thus $\phi\mapsto (\phi\cdot D)$ is a 
locally constant function on  $\arco(Z\subset X)$.

Note also that if $\phi\in \sharc\bigl(Z\subset Y\bigr)$,
then $(\phi\cdot D)>0$ iff $\phi(0)\in D$. 

 Let now
$D:=\sum_{i\in I} D_i\subset Y$ be a reducible divisor.
Fix natural numbers $\{m_i: i\in I\}$ and set
$$
A( m_i:i\in I):=\bigl\{\phi: 
(\phi\cdot D_i)=m_i: i\in I\bigr\}\subset \arco\bigl(D\subset Y\bigr).
$$
Since each $(\phi\cdot D_i) $ is locally constant on
$\arco\bigl(D\subset Y\bigr) $, every $A( m_i:i\in I) $
is a union of some connected components of $\arco\bigl(D\subset Y\bigr) $
and $\arco\bigl(D\subset Y\bigr) $ is the disjoint union of all the
$A( m_i:i\in I) $ for all $(m_1,\dots)\in \n^{|I|}$. 

Furthermore, if $\phi\in \sharc\bigl(D\subset Y\bigr)$,
then  $$
\phi\in A( m_i:i\in I)\qtq{iff}
\phi(0)\in \bigcap_{i: m_i>0}D_i.
$$
Let   $W\subset \bigcap_{i: m_i>0}D_i$ be
a connected component and set
$$
SA(W, m_i:i\in I):=\bigl\{\phi: \phi(0)\in W\mbox{ and } 
(\phi\cdot D_i)=m_i: i\in I\bigr\}\subset \sharc\bigl(D\subset Y\bigr).
$$
By the above considerations,  $SA(W, m_i:i\in I) $ 
is a union of connected components of $\sharc\bigl(D\subset Y\bigr) $.
\end{say}

Usually the  $SA(W, m_i:i\in I) $  are disconnected, but
they are connected and irreducible in the following important special case
which gives our basic supply of
connected arc spaces.

\begin{prop}\label{snc.short.arcs.lem}
 Let $Y$ be a complex manifold and 
$D:=\sum_{i\in I} D_i\subset Y$ a simple normal crossing divisor.
Then 
\begin{enumerate}
\item the above $SA(W, m_i:i\in I) $ are connected and
\item  they give all the  connected 
components of $\sharc\bigl(D\subset Y\bigr) $.
\end{enumerate}
\end{prop}

Proof. Sending an arc to its center gives  a continuous map
$SA(W, m_i:i\in I)\to W$. Since $W$ is connected, it remains to show that
each point $p\in W$ has an open neighborhood whose preimage in 
$SA(W, m_i:i\in I) $ is connected.
As we noted in (\ref{conn.comp.say}), this is a local question near $W$,
hence it is enough to prove the lemma when
$Y=\dd^n$, $D_i=(x_i=0)$ and $W=\cap_{i:m_i>0} (x_i=0)$.
This follows from the product formula (\ref{snc.short.arcs.say}.4). \qed

\begin{rem} If $X$ is an algebraic variety then each $SA(W, m_i:i\in I) $
contains a strongly irreducible dense open subset.
This can be seen as follows. 
After replacing $X$ with a  dense,  open subset
we may assume that there is a finite morphism
$\pi:X\to \c^n$ such that $D_i=\pi^{-1}(x_i=0)$. 
It is easy to see that our arc spaces are irreducible 
if $X=\c^n$ and the $D_i$ are hyperplanes.
The preimages of the resulting finite dimensional complex subspaces
give the necessary finite dimensional complex subspaces of
 $SA(W, m_i:i\in I) $.
\end{rem}

\section{Quotient singularities}\label{sec.quot}

Let $G\subset  \GL(n,\c)$ be a finite subgroup.
Then $G$ acts on $\c^n$ and  $X=\c^n/G$ is called a quotient singularity.
The quotient $X$ does not determine $G$ uniquely but there is a
smallest choice given by $G=\pi_1\bigl(X\setminus \sing X\bigr)$. 
This $G$ is characterized by the property that it contains no
pseudo-reflections, that is, elements whose fixed point set is a hyperplane.

We also need to work with quotients of simple normal crossing pairs,
thus we consider groups $G$ that also stabilize
 some of the coordinate hyperplanes.
This forces part of $G$ to be diagonal.

\begin{say}[Quotient singularities]\label{quot.sing.say}
Write  
$$
\c^n=\c^r_{x_1,\dots, x_r}\times \c^s_{y_1,\dots, y_s}.
$$
Set $D_i=(x_i=0)\subset \c^n$ and $D=D_1+\cdots+D_r$.
Let $G\subset (\c^*)^r\times \GL(s,\c)$ be a finite subgroup without 
pseudo-reflections; then each of the $D_i$ is  $G$-invariant.
Let $F\subset \c^n$ be the union of the fixed point sets of the 
non-identity elements. Thus $F$ is a union of linear subspaces of codimension
$\geq 2$ and $G$ acts freely on $\c^n\setminus F$.

 $(0\in X):=(0\in \c^n)/G$ is a quotient singularity;
let $\pi:\c^n\to X$ be the quotient map. Note that  $\sing X=\pi(F)$.
Set $E_i:=\pi(D_i)$ and $E=E_1+\cdots+E_r$.
For uniformity of notation, we set $E:=\{0\}$ if $r=0$.

We have an exact sequence
$$
0\to \z^r\to \pi_1\bigl(X\setminus(E\cup \sing X)\bigr)
\stackrel{\tau}{\to} G\to 1
\eqno{(\ref{quot.sing.say}.1)}
$$
and the $G$-action on $\z^r$ is trivial.
Thus, for every $g\in G$, the preimage $\tau^{-1}\langle g\rangle$ 
is an Abelian group,
in fact isomorphic to $\z^r$. This shows that if
$\gamma_1, \gamma_2\in \pi_1\bigl(X\setminus(E\cup \sing X)\bigr)$
are conjugate then 
$\tau(\gamma_1), \tau(\gamma_2)\in G$ are
conjugate and if $\tau(\gamma_1)=\tau(\gamma_2)$
then $\gamma_1, \gamma_2$ are conjugate iff $\gamma_1= \gamma_2$.

With the above notation, let
 $$
\sharc^*(E\subset X)\subset \sharc(E\subset  X)
\eqno{(\ref{quot.sing.say}.2)}
$$
 denote the  
set of those arcs $\phi:\bdd\to X$ such that
$\supp \phi^{-1}(\sing X)=\{0\}$. If $(0\in X)$ is an isolated
singularity, which is the main case that we are
interested in,  then $\sharc^*(E\subset X)= \sharc(E\subset  X)$.
\end{say}

\begin{thm} \label{quot.sing.arcs.cor}
Let $\bigl(0\in (X,E)\bigr)=\bigl(0\in (\c^n,D)\bigr)/G$ be a
 quotient singularity as above.
Then the irreducible components of $\sharc^*(E\subset X)$ 
 are in one-to-one correspondence
with the conjugacy classes of $\pi_1\bigl(X\setminus(E\cup \sing X)\bigr)$. 
\end{thm}

It is interesting to connect the above result
with the McKay correspondence; see \cite{reid-mckay} for a survey 
and for further references.

\begin{say}[Short arcs and the McKay correspondence]\label{mckay.say}
Let $\bigl(0\in \c^n/G\bigr)$ be an isolated quotient singularity
and $p:Y\to \bigl(0\in \c^n/G\bigr)$  a resolution of singularities.
For any short arc $\phi:\bdd\to \bigl(0\in \c^n/G\bigr)$,
let $p^{-1}\circ \phi: \bdd\to Y$ be its lift.
As $\phi$ varies in a dense open subset of an 
irreducible component of $\sharc(0\in \c^n/G)$, the centers
of the lifts $p^{-1}\circ \phi$ sweep out a dense subset of an
irreducible subvariety  of $Y$. 
By Theorem \ref{quot.sing.arcs.cor}, the
irreducible components of $\sharc(0\subset \c^n/G)$
correspond to the conjugacy classes of $G$, hence 
we get a natural  map
$$
\bigl\{\mbox{conjugacy classes of $G$}\bigr\}\to
\bigl\{\mbox{subvarieties  of $Y$}\bigr\}.
$$
(The trivial conjugacy class corresponds to arcs that can be deformed
away from the origin. The centers of such deformations sweep
out a dense subset of $Y$. Thus it is natural to 
say that $Y$ itself corresponds to the trivial conjugacy class.)

Assume next that  $G\subset \SL(2,\c)$ and  let
$p:Y\to \bigl(0\in \c^2/G\bigr)$
be the minimal resolution. Since the Nash conjecture holds for surfaces,
each irreducible exceptional curve of $p$ corresponds to
an irreducible component of the space of formal arcs 
$\farc\bigl(0\in \c^2/G\bigr)$ and each 
irreducible component of the space of formal arcs also gives
an irreducible component of the space of short arcs.

Thus short arcs give a concrete realization of the 
McKay correspondence between
nontrivial conjugacy classes of $G$ and
exceptional curves of the minimal resolution.

The McKay correspondence is not fully established in
higher dimensions. Our results suggest a connection
between the McKay correspondence and
the Nash conjecture for 
higher dimensional quotient singularities.
This topic should be explored further.
\end{say}

Our proof of Theorem \ref{quot.sing.arcs.cor}  is not very illuminating.
We write down a complete list of the
irreducible components of $\sharc^*(E\subset X)$  and then
observe the one-to-one correspondence
with the conjugacy classes.

\begin{say}[Arcs on quotient singularities]\label{arcs.quot.sing.say}
Let $\phi:\bdd\to X$ be an arc in $\sharc^*(E\subset X)$.
Set $D(\phi)=\bdd\times_X\c^n$ with normalization $\bar D(\phi)$.
Then $\bar D(\phi)$ is a disjoint union of several disks
$\coprod_i \bdd(\phi)_i$ and $\phi$ lifts to
$\phi_i:\bdd(\phi)_i\to \c^n$.

The restriction of $\pi$ gives 
 $\pi_i:\bdd(\phi)_i\to \bdd$; these are
Galois covers with cyclic Galois group $C_i\subset G$.
The different $C_i$ are conjugates of each other.
Set $m:=|C_i|$. Fixing a coordinate $s_i$ on $\bdd(\phi)_i$
such that $s_i^m=\pi_i^*t$ determines a generator $g_i\in C_i$.

Thinking of $g_i$ as a representation $\rho_i:\z/m\to \GL(r+s, \c)$,
$\phi_i:\bdd(\phi)_i\to \c^n$ is $\rho_i$-equivariant.

Let  $\sharc^*(D\subset \c^n)\subset \sharc(D\subset  \c^n)$ denote the  
set of those arcs $\phi:\bdd\to \c^n$ such that
$\supp \phi^{-1}(F)=\{0\}$. Then $\phi_i\in \sharc^*(D\subset \c^n)$.

The construction can be carried out in any connected family of arcs,
except that we may get a monodromy action of the fundamental group
of the base on the set of the disks $\bigl\{\bdd(\phi)_i\bigr\}$.
 This can be eliminated after
a finite covering.

Conversely, given $g\in G$, let $m$ be the order of $g$ and 
$\sharc^*(D\subset \c^n)^{\rho(g)}\subset \sharc^*(D\subset \c^n)$
 the set of $\rho(g)$-equivariant arcs.
Thus we have  established the following description of
$\sharc^*(E\subset X)$.

\begin{enumerate}
\item Pick $g\in G$.
Set $m=\ord(g)$, write $g$ as
$$
g=\bigl(\epsilon^{a_1}, \dots, \epsilon^{a_r}, g_s\bigr)
\qtq{where} g_s\in \GL(s,\c)
$$
and let $\rho(g):\z/m\to\GL(r+s, \c)$  denote the
corresponding representation.
Note that $(a_1, \dots, a_r)$ are unchanged by conjugation.
\item Pick $m_1,\dots, m_r\in \n$, not all zero, such that
 $m_i\equiv a_i \mod m$. 
\item By (\ref{snc.short.arcs.lem}) the $m_i$ determine 
$SA(m_1,\dots, m_r)\subset \sharc(D\subset \c^n)$ which is a
connected and irreducible component.
\item  Let $SA(m_1,\dots, m_r)^{\rho(g)}\subset SA(m_1,\dots, m_r)$
denote the subset of $\rho(g)$-equivariant arcs.
By (\ref{snc.short.arcs.say}.3)
and (\ref{arc.disk.say}),  $SA(m_1,\dots, m_r)^{\rho(g)}$
 is connected and irreducible.
\item We saw in (\ref{shrac*.prop}) that 
$$
SA^*(m_1,\dots, m_r)^{\rho(g)}:=SA(m_1,\dots, m_r)^{\rho(g)}
\cap \sharc^*(D\subset \c^n)
$$
is also connected and irreducible.
\item These arcs descend to arcs on $X$, giving
$$
 SA^*(m_1,\dots, m_r)^{\rho(g)}\to  \sharc^*(E\subset X).
$$
\item The image is a connected component
 of $\sharc^*(E\subset X)$.
\item Conjugate elements give the same 
connected component and every connected component
of $\sharc^*(E\subset X)$ is obtained this way.
\end{enumerate}
\end{say}

Next we  make this more explicit in two cases that are of 
special interest.

Assume first that $(0\in X)=(0\in \c^n)/G$ is an
 isolated quotient singularity and there are no divisors.
Thus $r=0$. Given $g\in G$ of order $m$,
$A^{\rho(g)}\subset \sharc(0\in \c^n)$ is the set of
$\rho(g)$-equivariant arcs. This $A^{\rho(g)}$ is connected and irreducible
by (\ref{snc.short.arcs.lem}).
Thus we obtain a proof of the following restatement of 
Theorem \ref{sharc.quot.thm.1}.

\begin{cor} \label{quot.sing.arcs.cor.cor}
Let $(0\in X)=(0\in \c^n)/G$ be an
 isolated quotient singularity where $G\subset \GL(n,\c)$ acts
freely outside the origin.
Then the irreducible components of $\sharc(0\in X)$ 
 are in one-to-one correspondence
with the conjugacy classes of $G$. \qed
\end{cor}

We also need to understand cyclic quotients of
$\bigl((x=0)\subset \c^2\bigr)$.
Consider the pair
$$
(E\subset X):=\bigl((x=0)\subset \c^2\bigr)/\tfrac1{m}(q,1).
$$
(The notation represents the quotient of $\bigl((x=0)\subset \c^2 \bigr) $
by the group action
$(x,y)\mapsto (\epsilon^q x, \epsilon y)$ where
$\epsilon=e^{2\pi i/m}$.)
Then $\pi_1(X\setminus E)$ is an extension of $\z/m$ by
$\z\cong \pi_1(\c^2\setminus D)$ where $D=(x=0)$. 
The universal cover of $\c^2\setminus D $ is
$$
\rho:\c^2_{uv}\to \c^2_{xy}\setminus D
\qtq{where} \rho(u,v)=\bigl(u, e^{2\pi i v}\bigr).
$$
The group of deck transformations is generated by
$(u,v)\mapsto (u, v+1)$.
We see that $\c^2_{uv}$ is also the universal cover of
$X\setminus E $ and the 
group of deck transformations is generated by
$$
(u,v)\mapsto \bigl(e^{2\pi i q/m} u, v+\tfrac1{m}\bigr).
$$
Given $a\in \{0,\dots, m-1\}$ (which we think of as an element of
$\z/m$) let $m_1$ be a positive integer 
congruent to $a\mod m$ and 
 $c$  the smallest nonnegative integer such that 
$a\equiv cq\mod m$. Then the arc
$$
\phi(m_1):\bdd\to \c^2\qtq{given by}  t\mapsto \bigl(t^{m_1}, t^{c}\bigr)
$$
is a typical element of  $SA(m_1)$ in (\ref{arcs.quot.sing.say}.3). 
The intersection number of $\phi(m_1) $ with  $D=(x=0)$
is $m_1$. Thus $\phi(m_1)$ descends to an arc in $X$ whose
intersection number with $E$ is $m_1/m$. Thus we obtain the following.

\begin{cor} \label{2dim.dlt.arcs.cor}
Let $(E\subset X):=\bigl((x=0)\subset \c^2\bigr)/\tfrac1{m}(q,1) $.
The intersection number establishes a bijection 
$$
\pi_0\bigl(\sharc(E\subset X)\bigr)\leftrightarrow \tfrac1{m}\z_{>0}.
$$
If $a\in \tfrac1{m}\z_{>0}$ is not an integer, then
the corresponding arcs are all centered at the origin.
If $a\in \z_{>0}$ then a corresponding general arc 
has center on $E\setminus \{0\}$. \qed
\end{cor}

\section{Minimal dlt modifications}\label{sec.dlt}

As outlined in Paragraph \ref{method.sharc.surf.thm.2},
the first step toward proving Theorem \ref{sharc.surf.thm.2}
is the construction and study of minimal dlt modifications
of surfaces.

\begin{defn}[Minimal dlt modification]\label{dlt.surf.defn}
Let $Y$ be a normal surface and $E\subset Y$ a reduced curve.
The pair $(Y,E)$ is {\it divisorial log terminal}
(abbreviated as {\it dlt}) 
 if everywhere it has one of the following
normal forms (analytically locally).
\begin{enumerate}
\item (Normal crossing points)
$$
\bigl((xy=0)\subset \c^2 \bigr)\qtq{or} \bigl((x=0)\subset \c^2 \bigr)
\qtq{or} \bigl(\emptyset\subset \c^2 \bigr).
$$
\item (Cyclic quotients)
$$
\bigl((x=0)\subset \c^2 \bigr)/\tfrac1{m}(q,1)\qtq{where} (q,m)=1.
$$
\item (Quotients) $\bigl(\emptyset \subset \c^2\bigr)/G$
where $G\subset \GL(2,\c)$ is a finite subgroup 
acting freely outside the origin.
\end{enumerate}
(For a more conceptual definition, see \cite[2.37]{km-book}.) 

Let $(0\in X)$ be a normal surface singularity
and $f:Y\to X$ a proper birational morphism
that is an isomorphism outside $0$. Set $E:=\supp f^{-1}(0)$.
We say that  $f:(E\subset Y)\to (0\in X)$ is a 
{\it  dlt modification} if
\begin{enumerate}\setcounter{enumi}{3}
\item  $(Y,E) $ is dlt and
\item  $\bigl((K_Y+E)\cdot E_i\bigr)\geq 0$  for every 
irreducible curve $E_i\subset E$.
\end{enumerate}
A dlt modification is called {\it minimal} if, in addition,
\begin{enumerate}\setcounter{enumi}{5}
\item there are no exceptional curves
 $E_i\subset E$ such that $Y$ is smooth along $E_i$,
$\p^1\cong E_i$,  $(E_i^2)=-1$,  $\bigl((K_Y+E)\cdot E_i\bigr)= 0$ and 
$E_i$ intersects the rest of $E$ 
in  2   points.
\end{enumerate}

\end{defn}

\begin{say}[Construction of minimal dlt modifications]\label{min.dlt.mod.defn}
Let $(0\in X)$ be a normal surface singularity.

If $X$ is a quotient singularity, then it is its own
minimal dlt modification. 

If $X$ is not a quotient singularity, then its 
minimal dlt modification can be constructed as follows.

Let $g:(E\subset Y)\to (0\in X)$  be a log resolution,
that is, $Y$ is smooth and $E:=\supp g^{-1}(0)$ is a
curve with nodes only. If $\p^1\cong E_i\subset E$ is a $-1$-curve
that intersects the rest of $E$ 
in at most 2  distinct points, 
 then we can contract $E_i$ to get another  
log resolution. After all such contractions we get
the {\it minimal log resolution}
$$
g^m:(E^m\subset Y^m)\to (0\in X).
$$
Next we contract maximal tails of rational chains. 
A sequence of curves  $E_1,\dots, E_m$ is called a
rational chain if these are smooth rational curves
and 
$$
(E_i\cdot E_j)=
\left\{
\begin{array}{l}
0 \qtq{if} |i-j|>1\qtq{and} \\
1 \qtq{if} |i-j|=1.
\end{array}
\right.
$$
The sequence is called a tail if it intersects the rest of $E$ in a single 
point of $E_m$.  The maximal  tails of rational chains are
disjoint from each other. If $(0\in X)$ is not a  quotient  singularity
then contracting each of the  maximal  tails of rational chains  gives
the  {\it minimal dlt modification}
$$
g^{\rm dlt}:(E^{\rm dlt}\subset X^{\rm dlt})\to (0\in X).
$$

Note that if $(0\in X)$ is a cyclic quotient  singularity 
then the latter recipe gives
the correct answer: $X^{\rm dlt}=X$.  
If $(0\in X)$ is a non-cyclic quotient  singularity the above recipe gives
 $X'\to X$ with a unique exceptional curve $E'$ and
3 singular points on $X'$. This is not 
a dlt modification since 
$\bigl(\bigl(K_{X'}+E'\bigr)\cdot E'\bigr)<0$.
The correct  minimal dlt modification is again $X^{\rm dlt}=X$.
\end{say}

Following the method of (\ref{snc.short.arcs.lem}),
 we get the following description
of  the arc space of dlt pairs. The only difference is that we have to
understand arcs centered at the singular points;
these were described in Corollary \ref{2dim.dlt.arcs.cor}.

\begin{say}[Short arcs on dlt pairs]\label{sharc.dlt.pairs.say}
Let $(E\subset Y)$ be a dlt pair. Then the following is a list
of the connected components of $\sharc(E\subset Y)$.
\begin{enumerate}
\item Let $E_i$ be an irreducible component of $E$ and $m_i\in \z_{>0}$.
We get $SA(E_i, m_i)$, consisting of arcs whose intersection number with
$E_i$ is $m_i$ (and with every other $E_j$ is $0$).
\item Let 
$E'_p, E''_p$ be the two local branches of $E$ at  a singular point $p$
  and $m'_p, m''_p\in \z_{>0}$.
We get $SA(p, m'_p, m''_p)$, consisting of arcs whose 
center is at $p$ and whose 
intersection numbers with
the two local branches  are $m'_p, m''_p$.
\item Let $q$ be a singular point of $X$ of the form
$$
\bigl((x=0)\subset \c^2 \bigr)/\tfrac1{m}(b,1)
\qtq{and}
m_q\in \tfrac1{m}\z_{>0}\setminus \z_{>0}.
$$
We get $SA(q, m_q)$, consisting of arcs whose 
center is at $q$ and whose 
intersection number with $E$ is $m_q$.
\end{enumerate}
In the above cases it is easy to work out the homotopy types of the
different connected components of $\sharc(E\subset Y)$.
\begin{enumerate}
\item[(1')] Let $E_i^0\subset E_i$ be the complement of the 
set of singular points of $E$. Then $SA(E_i, m_i)$ is
homotopy equivalent to an $\s^1$-bundle over $E_i^0$.
If $E_i$ has genus $g$ and contains $r>0$ singular points
then  $SA(E_i, m_i)$ is
homotopy equivalent to
$\s^1\times \bigvee_1^{2g+r-1}\s^1$. The exceptional case is when
 $E=E_i$ consists of a single smooth curve, thus $r=0$.
Then we get an $\s^1$-bundle over a compact Riemann surface (without boundary)
whose Chern class is $m_i(-E_i^2)$. 
\item[(2')] Taking the angle component of the
leading coefficients of the Taylor series
of an arc shows that $SA(p, m'_p, m''_p)$ is
homotopy equivalent to  $\s^1\times \s^1$.
\item[(3')] Working with equivariant  Taylor series on the universal cover
shows that $SA(q, m_q)$  is
homotopy equivalent to  $\s^1$.
\end{enumerate}
These cases are topologically distinct, except that $SA(E_i, m_i)$
is homotopy equivalent to  $\s^1\times \s^1$ if $g=0$ and $r=2$. 
\end{say}

Our main result says that (\ref{sharc.dlt.pairs.say}.1--3) gives a
complete list of the connected components of $\sharc\bigl(0\in X\bigr)$
for normal surface  singularies.

\begin{thm} \label{dlt.nonquot.thm}
 Let $(0\in X)$ be a normal surface  singularity  
with  minimal dlt modification
$g^{\rm dlt}:\bigl(E^{\rm dlt}\subset X^{\rm dlt}\bigr)\to (0\in X) $.
 Then composing an arc with
$g^{\rm dlt} $ gives a bijection
$$
g^{\rm dlt}_*: 
\pi_0\Bigl(\sharc\bigl(E^{\rm dlt}\subset X^{\rm dlt}\bigr)\Bigr)
\ \leftrightarrow\
\pi_0\Bigl(\sharc\bigl(0\in X\bigr)\Bigr).
\eqno{(\ref{dlt.nonquot.thm}.1)}
$$
\end{thm}

Proof. If $(0\in X)$ is 
a   quotient  singularity then $X^{\rm dlt} =X$ and there is nothing to prove.

Since $f_*$ gives a continuous bijection
$$
f_*:\sharc\bigl(E^{\rm dlt}\subset Y^{\rm dlt}\bigr)\to \sharc(0\in X),
$$
the induced map on $\pi_0$ is surjective and it remains to prove that
it is injective.  Using the winding number maps (\ref{conn.comp.say}.2--3),
(\ref{dlt.nonquot.thm}.1) sits in a diagram
$$
\begin{array}{ccc}
\pi_0\Bigl(\sharc\bigl(E^{\rm dlt}\subset X^{\rm dlt}\bigr)\Bigr) 
& \stackrel{w}{\to}&
\pi_1\Bigl(\link\bigl(E^{\rm dlt}\subset X^{\rm dlt}\bigr)\Bigr)/(\mbox{conj})\\
g^{\rm dlt}_*\downarrow && \downarrow \cong\\
\pi_0\Bigl(\sharc\bigl(0\in X\bigr)\Bigr) & \stackrel{w}{\to} &
\pi_1\bigl(\link(0\in X)\bigr)/(\mbox{conj}).
\end{array}
$$
Thus it is sufficient to prove that  the winding number map
$$
w:\pi_0\Bigl(\sharc\bigl(E^{\rm dlt}\subset X^{\rm dlt}\bigr)\Bigr)\to
\pi_1\Bigl(\link\bigl(E^{\rm dlt}\subset X^{\rm dlt}\bigr)\Bigr)/(\mbox{conj})
$$
is an injection. The latter is proved in Section \ref{sec.fundgr}. \qed 

\begin{rem} \label{dlt.nonquot.rem}
 It is possible that the map
$$
g^{\rm dlt}_*: 
\sharc\bigl(E^{\rm dlt}\subset Y^{\rm dlt}\bigr)
\to
\sharc\bigl(0\in X\bigr)
\eqno{(\ref{dlt.nonquot.rem}.1)}
$$
is a homeomorphism but we do not know how to prove this.
\end{rem}

\section{Fundamental groups of normal surface 
singularity links}\label{sec.fundgr}

For general introductions to the topics of 3--manifold topology that we  use,
see \cite{seif, MR0415619, seif-book, MR565450, scott}.
Connections with links are treated in detail in  \cite{NP}.

In this section we denote by $L$ the oriented link of a 
normal surface singularity
$(0\in X)$, and by $\pi_1$ its
fundamental group. Recall that $L$ is connected.
\begin{thm}\label{pi1.prop} The following properties hold.
\begin{enumerate}
\item  Singularity links are irreducible 3--manifolds (that is,
every embedded $\s^2$ in $L$ bounds a 3--ball).
\item $\pi_1$ is finite iff $X$ is a quotient singularity.
\item If $\pi_1$ in infinite
then  $L$ is the Eilenberg--Mac Lane space of $\pi_1$, that is $L=K(\pi_1,1)$.
Moreover, in this case, $\pi_1$ is torsion free.
\end{enumerate}
\end{thm}
Proof. Part (1) follows e.g. from  \cite[Thm.1]{NP}. If $\pi_1$ is finite then the
universal cover of $L$ is $\s^3$, the link of $\c^2$ at the origin,
and the action of the finite $\pi_1$ can be linearized. For (3) see e.g.
\cite[p.32]{3man}. The argument is the following.
By the Sphere Theorem, if $N$ is an orientable 3--manifold with 
$\pi_2(N)\not=0$,
then $N$ contains an embedded 2--sphere which is homotopically non--trivial.
Hence, since $L$ is irreducible, we must have $\pi_2(L)=0$. Since $\pi_1$ in infinite,
the 3--manifold $\widetilde{L}$, the universal cover of $L$, has $H_3(\widetilde{L})=0$.
Hence, by the Hurewicz theorem 
$\pi_i(L)=\pi_i(\widetilde{L})=H_i(\widetilde{L})=0$ for any $i>1$.
Next, any finite dimensional Eilenberg--Mac Lane space $L=K(\pi_1,1)$ has a torsion free fundamental group.
\qed

\vspace{2mm}

Part (3) can also be proved using the JSJ 
(=Jaco--Shalen--Johannson)
decomposition of $L$ (see below), and using Bass--Serre theory
(for the  $K(\pi_1,1)$ property) and the Torsion Theorem for free products with amalgamations and
HNN (=Higman--Neumann--Neumann)
extensions (for the torsion free property). Similar arguments will be used in the next
paragraphs for conjugacy properties.

The case of finite $\pi_1$ is completely classified \cite{bri-rat,DuVal},
they are the quotient singularities.  The
relevant statement regarding $\sharc$ is given in Section \ref{sec.quot}.

In the sequel we will assume that $\pi_1$ is an infinite group.

The fundamental group $\pi_1$ has a presentation in terms of the plumbing realization
of the link via a dual resolution graph. This presentation usually is rather 
involved due to the cycles in the graph.
Hence, for the first conceptual conjugacy statements,
we will use the  structural decomposition of $\pi_1$ as the fundamental group of
a graph of groups induced by the JSJ decomposition. Then, for each independent
component of the JSJ  decomposition (whose graphs have no cycles anymore)
we will consider the corresponding concrete  plumbing presentation.

\begin{say}[Plumbing graphs]\label{say:pi1_pl}
Let us fix the minimal log resolution of
$(0\in X)$ as in (\ref{min.dlt.mod.defn}). Then the link appears as
the boundary of a small tubular neighbourhood of the exceptional curve, or,
as a plumbed 3--manifold associated with the
dual resolution graph $\Gamma$ of the irreducible exceptional components.
The  set of vertices of this graph will be denoted by $\sV$, the set of
edges by $\sE$.
The {\it nodes} are  those vertices $v$ which either have 
genus $g_v>0$ or valency $\geq 3$.

The plumbing description of the link provides a presentation of $\pi_1$.
For rational homology sphere links (when the exceptional curve has trivial  fundamental group,
and $\pi_1$ is generated by the loops around the irreducible exceptional components)
this was determined by Mumford \cite{mumf-top},  for the general case see
\cite{catanese-MR2238372, ausina}. We will use this concrete presentation only for the
  Seifert pieces, cf.\ (\ref{say:pi1j}).

Each edge of the plumbing provides naturally an embedded  torus in $L$.
The tori of the minimal JSJ decomposition are such tori: from each maximal chain with node ends
(which can be the very same node in the case of a `loop') we choose
one edge, their collection serves as the JSJ--tori for the minimal decomposition.
(For a fixed such  chain, the choice of the edge will not alter the torus, but only a canonical basis
in its homology $\z^2$.)
In particular, the pieces $\{L_j\}_j$ of the decomposition are indexed by nodes, and each piece has a Seifert structure whose central vertex is the corresponding node.

In fact, there is an exception to this description, when $\Gamma$ is a cyclic graph with  all
genus decorations zero (the case of cusps): one has no nodes, exactly one JSJ torus provided by one of the edges arbitrarily chosen,  and one Seifert piece.

In the plumbing representation of the pieces $L_j$ we extend the above notation
of the plumbing graphs:
 the boundary components of $L_j$ will be denoted by {\it arrowhead} vertices
and their set is denoted by  $\sA$. The set of the other, 
 non--arrowhead vertices,  is denoted  by $\sW$. Hence $\sV=\sW\sqcup \sA$. 
If an arrow $a$ is supported by
 the non--arrowhead $v$, then the corresponding boundary component is the boundary of the tubular
 neighborhood of a generic $\s^1$--fiver $S_a$ over the surface (irreducible exceptional curve)
 indexed by $v$ in the plumbing construction.
The plumbing graph of the collection of $\{L_j\}$ is obtained from the graph of $L$ by
replacing each chosen edge of the chains by two arrowheads (supported by the two endpoints of the edge).

In the remaining part of this section we assume that $\Gamma$ has at least one node.
The case of cyclic graphs without nodes will be treated in Section \ref{sec.cusps}.
\end{say}
\begin{say}[The non--Seifert case. The presentation of $\pi_1$ via graph of groups]\label{say:pi1}
Assume first that $L$ is not a Seifert 3--manifold and
 $\Gamma$ has at least one node. Then, being irreducible, and having at least one
incompressible  JSJ torus, it is Haken.
 The JSJ decomposition of $L$ provides a splitting of $\pi_1$
as a fundamental group of a graph of groups \cite{serre-trees,Pr}.
More precisely, the JSJ tori form a collection of disjointly embedded incompressible tori
 $\{T_\tau\}_\tau$ (where incompressibility means that $\pi_1(T_\tau)\to \pi_1(L)$ is injective),
 such that the components $\{L_j\}_{j}$
 of $L\setminus \cup_{\tau}T_\tau$ are Seifert fibred.
 Each $T_\tau$ is two--sided in $L$, let
  $T_\tau^+$ and $T_\tau^-$
 be the two boundaries of a small tubular neighbourhood of $T_\tau$ in $L$.

The underlying graph ${\mathfrak G}$
of the graph of groups has one vertex $v_j$ for each piece $L_j$, we  will also
write $L(v_j)=L_j$. To each $T_\tau$ correspond two edges of ${\mathfrak G}$,
$e_\tau:=e(T_\tau)$ and $\overline{e_\tau}:=
e(T_\tau)^{-1}$, inverse to each other. If $L_k$ and $L_l$ are such that $T_\tau^-\subset L_k$ and $T_\tau^+\subset
L_l$, then we say that $e(T_\tau)$ has origin $v(L_k)$ and extremity $v(L_l)$.

The graph of groups is obtained by assigning to each vertex $v_j$ the {\it vertex group}
$G(v_j)=\pi_1(L(v_j))$, and to each edge $e_\tau$ the {\it edge group} $G(e_\tau)=\pi_1(T_\tau)
\simeq\z\oplus \z$.
For each edge $e_\tau$ with origin $v_k$ and extremity $v_l$, one defines the monomorphisms
$\phi_\tau^-:G(e_\tau)\hookrightarrow G(v_k)$ and
$\phi_\tau^+:G(e_\tau)\hookrightarrow G(v_l)$.
(Strictly speaking, this depends on the choice of base points $*_j$ in $L_j$, $*_\tau^\pm\in T_\tau^\pm$,
and paths from $*_\tau^-$ to $*_k$ and from $*_\tau^+$ to $*_l$ respectively.)
Their images in the corresponding vertex groups will be denoted by
$G(e_\tau)^-$ and $G(e_\tau)^+$, and we set $\phi_{e_\tau}=\phi_\tau^+\circ (\phi_\tau^-)^{-1}:G(e_\tau)^-\to
G(e_\tau)^+$. For the inverse edges  $G(\bar{e})=G(e)$, $G(\bar{e})^{\pm}=G(e)^{\mp}$,
$\phi_{\bar{e}}=\phi_e^{-1}$.

These data determine $\pi_1$ as the fundamental group of the graph of groups,
once a maximal tree ${\mathfrak T}$ in ${\mathfrak G}$ is chosen. An edge of ${\mathfrak G}$
is called ${\mathfrak T}$--separating if it belongs to
${\mathfrak T}$, and ${\mathfrak T}$--non--separating otherwise.
Then, if for each vertex $v$,  $\langle S_v|R_v\rangle$
is the presentations of $G(v)$, then the presentation of $\pi_1$ is
the following (cf.\ \cite[5.1]{serre-trees} or
\cite[page 179]{Pr}):

 \vspace{2mm}

 Generators: $(\cup_{j=1}^n S_{v_j})\cup \{t_e| \, e\, \mbox{edge of ${\mathfrak G}$}\}$.

 Relations: \begin{equation*}\begin{split}
 (\cup_{j=1}^n R_{v_j})& \cup
 \{ \mbox{for all edge $e$ of ${\mathfrak G}$ and generator  $c$ of $G(e)^-: \ t_e\phi_e(c)t_e^{-1}=c$}\}\\
 & \cup \{ \mbox{for all edge $e$ of ${\mathfrak G}$: \ $t_{\bar{e}}=t_e^{-1}$}\}\\
  & \cup \{ \mbox{for all ${\mathfrak T}$-separating edge $e$ of ${\mathfrak G}$: \
  $t_e=1$}\}.\end{split}\end{equation*}
\end{say}
The presentation $\langle S_v|R_v\rangle$ of each contributing group  $G(v)$ will be given  in
(\ref{say:pi1j}).

\begin{say}[Reduced and cyclically reduced forms]\label{say:words} \cite{serre-trees,Pr}
A {\it path} in the graph ${\mathfrak G}$, of length $m\geq 0$,  is a sequence $(v_{j_0},e_{\tau_0},v_{j_1},\ldots,
e_{\tau_{m-1}},v_{j_m})$, where the origin and the extremity of $e_{\tau_i}$ are $v_{j_i}$ and
 $v_{j_{i+1}}$ respectively. A path is called a {\it loop} if $j_0=j_m$.

A {\it word of type $\sC$} is a couple $(\sC,\mu)$, where $\sC$ is a loop in ${\mathfrak G}$, and
$\mu$ is a sequence $\mu=(\mu_0,\ldots, \mu_m)$, such that $\mu_i\in G(v_{j_i})$ for all $i$.

Once a base point (vertex) $v_*$ in ${\mathfrak G}$ is fixed, we might consider loops with  $v_{j_0}=v_{j_m}=v_*$.
Then any such word $(\sC,\mu)$ gives an element of the fundamental group of the graph of groups, denoted by
$|\sC,\mu|$. In the presentation (\ref{say:pi1}) it is the word $\mu_{j_0}\cdot t_{\tau_0}
\cdot \mu_{j_1}\cdots  t_{\tau_{m-1}}\cdot \mu_{j_m}$.
We say that $|\sC,\mu|$ is  the label of $(\sC,\mu)$.

 A word $(\sC,\mu)$ is called {\it reduced} if either $m=0$ and $\mu_{j_0}\not=1$, or $m>0$ and
 any time when $e_{\tau_{i-1}}=\overline{e_{\tau_i}}$
 one has $\mu_i\in G(v_{j_i})\setminus G(e_{\tau_i})^-$.
 Given any non-trivial element of the fundamental group, there exists a reduced form associated with it.

 A {\it cyclic conjugate} of the word $(\sC,\mu)$ is $(\sC',\mu')$, where
 $$\sC'=(v_{j_i},e_{\alpha_i},v_{j_{i+1}},\ldots, e_{\tau_{m-1}},v_{j_m}, e_{\tau_0},
 \ldots , e_{\tau_{i-1}},v_{j_i})$$ and
 $\mu'=(\mu_i,\mu_{i+1},\ldots,\mu_{m-1},\mu_m, \mu_0,\ldots, \mu_{i-1},1)$, for some $1\leq i\leq m$.

 A word is called {\it cyclically reduced}
  if all its cyclic conjugates are reduced, and if $\mu_m=1$ whenever $m>0$.
  Any non-trivial  conjugacy class of the fundamental group  can be represented by a cyclically reduced
  form whose label is an element of the class.

  We say that a word $(\sC,\mu)$ is {\it 0--reduced}, if $\mu_i=1$ for all indices $i$, except one.
  Note that a 0--reduced word has a cyclically reduced form of length 0.
  \end{say}
  \begin{say}
 The {\it Conjugacy Theorem in the  group of graph of groups}
 characterizes (cyclically reduced) words whose  labels are conjugate
 in the group; see e.g. \cite[Theorem 3.1]{Pr}.
 Since in our application we need the conjugacy properties  of special 0--reduced
  words, we state the theorem only for them.
Below $\sim$ denotes conjugacy.
\end{say}
\begin{thm}\label{th:conj} \cite{Pr}
Suppose that $(\sC,\mu)=(v_0,\mu_0)$
 and $(\sC',\mu')=(v_0',\mu_0')$ are two cyclically reduced forms, both of length zero, whose labels
 $\mu_0$ and $\mu_0'$ are conjugate in $\pi_1$. Then  
 \begin{enumerate}
 \item either $v_0=v_0'$, and the elements  $\mu_0$ and $\mu_0'$ are conjugate in $G(v_0)$, 
 \item or  there exists a path
 $(v_{j_0},e_{\tau_0},v_{j_1},\ldots, e_{\tau_{m-1}},v_{j_m})$
  of length $m>0$,  with
 $v_{j_0}=v_0$ and $v_{j_m}=v_0'$, and a sequence of elements $(c_0,\cdots , c_{m-1})$
 with $c_i\in G(e_{\tau_i})$, such that :
  \begin{equation*}\begin{split}
  \mu_0 & \ \sim \ \phi_{e_{\tau_0}}^-(c_0) \ \mbox{in} \ G(v_{j_0}),\\
  \phi_{e_{\tau_{i-1}}}^+(c_{i-1})& \ \sim \ \phi_{e_{\tau_i}}^-(c_i) \ \, \mbox{in} \ G(v_{j_i}), \
  1\leq i\leq m-1,\\
 \phi_{e_{\tau_{m-1}}}^+(c_{m-1})& \ \sim \   \mu_0' \ \ \ \mbox{in} \ G(v_{j_m}).
  \end{split}\end{equation*}
 \end{enumerate}
\end{thm}

\begin{say}[Arc--generators and m-arc-generators]\label{say:arcgens}
 Consider again the plumbing representation of $L$, cf.\  (\ref{say:pi1_pl}).
 For a vertex $v\in \sV$ of the plumbing graph $\Gamma$  let $\g_v$ be the
oriented $\s^1$-fiber associated with $v$ in the plumbing construction.
Note that 
$\g_v$ is well defined up only to conjugacy; nevertheless, 
if we fix a base point $*$ of $L$ and a
connecting path from $*$
to a point of this loop, we get an element of $\pi_1$.
If $u,v$ are connected by an edge in $\Gamma$,
then we can consider for both of them the same connecting path,
in which case  $\g_u$ and $\g_v$ commute in $\pi_1$. (Similarly one can be achieve that all the
$\g$--generators along a chain commute with each other; see also the relations between them in
the second paragraph of (\ref{say:pi1j}).)

Similarly, assume that the  JSJ
decomposition is fixed (that is, the separating edges of the plumbing graph are fixed).
 Then $\g_v$ for each
 $v\in\sV(L_j)$ determines a well-defined conjugacy class of
 $\pi_1(L_j)$. Moreover,
 for each arrowhead $a\in\sA(L_j)$, inherited from the edge $e_\tau$ such that
  $L_j\supset T_\tau^\circ$, $\circ =+$ or $-$, we can consider $\g_a$, the oriented loop
 around $S_a$. If $v(a)$ is the supporting vertex of $a$ in the graph of $L_j$, then by taking for
 the connecting path of $\g_a$ and $\g_{v(a)}$ the path connecting the basepoints $*_j$ and $*_\tau^\circ$
 already considered in (\ref{say:pi1}), we obtain that  $\g_a$ and $\g_{v(a)}$ generate $G(e_\tau)^\circ \simeq\z\oplus\z\subset \pi_1(L_j)$.

 In the sequel we abbreviate $\pi_1(L_j)$ by $\pi_{1,j}$ and $G(e_\tau)^\circ$ by $H_\tau\subset \pi_{1,j}$.

We consider the following conjugacy classes of $\pi_1$:
 for each $v\in \sV$ we take the class of  $\g_v^{m_v}$ for any fixed $m_v>0$, and
for each pair $v,u\in\sV$, $(u,v)\in\sE$,  the class $\g_u^{m_u}\g_v^{m_v}$ for  $m_u,\,m_v>0$
(with common connecting paths for $u,v$).
We call these classes {\it arc-generators}.
Their definition is motivated by the map $\sharc\to \pi_1/({\rm conj})$,
see (\ref{conn.comp.say}.2) and (\ref{sharc.dlt.pairs.say}).
It is convenient to extend this family and include also those powers
which have negative exponents: we call these
classes  {\it m-arc-generators} (`m' from `meromorphic').
That is,  they consists of classes of type
 $\g_v^{m_v}$, with $m_v\in\z$ and $\g_u^{m_u}\g_v^{m_v}$ with $m_u,\,m_v\in\z$ (for $(u,v)\in\sE$).
Thus we conclude the following.
\end{say}
\begin{lem}\label{lem:marc}
The conjugacy classes of m-arc-generators of $\pi_1$ can be represented by 0-reduced words, hence their
cyclically reduced forms have length zero.\qed
\end{lem}

In particular, for conjugate  arc-generators Theorem \ref{th:conj} can be applied.
Moreover, we will show that for  arc-generators necessarily case (1) happens, hence we have to
understand  conjugate arc-generators of some fixed
$\pi_{1,j}$ only.  For this
 we need the precise presentation of $\pi_{1,j}$.

\begin{say}[The presentation of $\pi_{1,j}$]\label{say:pi1j}
The plumbing graph of $L_j$ is star--shaped.
Let the number of  chains without
arrows at the end be $n\geq 0$, and the number of chains with
arrow at the end be $k\geq 1$.

Each chain without arrowhead at the end has a Seifert invariant
$(\alpha_i,\omega_i)$, $0<\omega_i<\alpha_i$,
${\rm gcd}(\omega_i,\alpha_i)=1$. If $[b_1,\ldots, b_s]$ is the Hirzebruch continued fraction
of $\alpha_i/\omega_i$, then the chain has $s$
vertices, all of them of genus 0, their Euler numbers are  $-b_1,\ldots , -b_s$ such that the
vertex decorated with $-b_1$ is connected with the node.
If $\g_1,\ldots, \g_s$ are the corresponding elements  in the fundamental group, then
one has the relations $\g_s^{b_s}=\g_{s-1}$, \ldots, $\g_{i}^{b_i}=\g_{i-1}\g_{i+1}$ (see e.g. \cite{mumf-top}).
In particular, all the generators $\g_i$ ($1<i\leq s$)
  can be eliminated: they appear as certain  powers of $\g_s$. Indeed,
let $\det[b_k,\ldots, b_l]$ be the numerator of the corresponding
continued fraction (or the determinant of the corresponding
subgraph), then $\g_i=\g_s^{\det[b_s,\ldots,b_{i+1}]}$.

Similar elimination happens for a chain with arrow at the end,
but now  the  generator of that
vertex survives which is connected with the node
 (independently whether it is arrowhead or non--arrowhead), 
see also the proof of Lemma \ref{lem:z2}.
(Note also that the Euler numbers along these arrowed chains become irrelevant
in the topology of $L_j$ --- with non-framed boundary ---, and in the presentation of $\pi_{1,j}$ as well.)

We rename the non-eliminated  generators: $h$ is the generator of the node, $\{g_i\}_{i=1}^n$ are
the generators of the end vertices of chains without arrow, and
$\{f_\tau\}_{\tau=1}^k$ are the generators of vertices next to the node
on the chains with arrows. Let $-b$ and $g$ be the  Euler number and  genus of the node.
Then, if $g\not=0$ we will have $2g$ more generators, $\{a_m,b_m\}_{m=1}^g$, and
$\pi_{1,j}$ has the following presentation \cite{mumf-top,Pr}.
\begin{equation*}
\langle h, \{g_i\}_{i},\{f_\tau\}_{\tau},\{ a_m,b_m\}_{m} \
|\ h \ \mbox{central},\ g_i^{\alpha_i}=h,\forall i, \
h^b = \tprod_m [a_m,b_m]\cdot\tprod_i g_i^{\omega_i}\cdot \tprod_\tau f_\tau\rangle.
\end{equation*}
Note that this group is independent of $b$ since
 we can eliminate  one of the generators
 $f_{\tau_0}$ and the last equation.
Then the generators $\{f_\tau\}_{\tau\not=\tau_0}$ and $\{ a_m,b_m\}_{m}$ become free generators, and
$\pi_{1,j}$ is a free product $G*F_{2g+k-1}$, where $F_r$ is the free group with $r$ letters, and
$$
G:=\langle h, \{g_i\}_{i=1}^n  \ |\  g_i^{\alpha_i}=h, \forall i \rangle.
\eqno{(\ref{say:pi1j}.1)}
$$
Since the central element
 $h$ has infinite order, cf.\ (\ref{pi1.prop}), it generates a normal subgroup
$\langle h\rangle \simeq\z$. Moreover,
$$
\overline{\pi_{1,j}}:=\pi_{1,j}/\langle h\rangle \simeq\z_{\alpha_1}*\cdots *\z_{\alpha_n}*F_{2g+k-1}.
\eqno{(\ref{say:pi1j}.2)}
$$
\end{say}
\begin{notation} In the sequel
for any $x\in \pi_{1,j}$ we denote by $\bar{x}$ its image in $\overline{\pi_{1,j}}$, and
$\ell:=\tprod_m[a_m,b_m]\cdot \tprod_ig_i^{\omega_i}$.
\end{notation}
\begin{lem}\label{lem:z2} Fix an arrow (or boundary torus) indexed by $\tau$ ($1\leq \tau\leq k$).
Then  $h$ and $f_\tau$ generate $H_\tau\simeq \z\oplus\z\subset \pi_{1,j}$. (Recall
$H_\tau=G(e_\tau)^\circ$, cf.\ (\ref{say:arcgens}).)

In particular, the m-arc-generators of $L_j$ in $\pi_{1,j}$ are the following: either elements of
$H_\tau$ ($1\leq \tau\leq k$), or integral powers of $g_i$ ($1\leq i\leq n$).
\end{lem}
Proof.
Let $-b_1,\ldots, -b_s$ be the Euler numbers of the vertices of the corresponding chain,
where the $s$-th vertex supports
the arrow.  Via similar identities as in (\ref{say:pi1j}), any $\g_v$ along this chain
equals $f_{\tau}^ah^b$, where the exponents are certain determinants; in particular,
 with the notation of (\ref{say:arcgens}):
$$\g_a=f_\tau^{\det[b_1,\ldots, b_s]}h^{- \det[b_2,\ldots, b_s]}, \ \
\g_{v(a)}=f_\tau^{\det[b_1,\ldots, b_{s-1}]}h^{- \det[b_2,\ldots, b_{s-1}]}.$$
Since this transformation  from $(f_\tau,h)$ to $(\g_a,\g_{v(a)})$  has
unimodular matrix,
and $H_\tau$ is generated by $\g_a$ and $\g_{v(a)}$, the lemma follows. \qed

\begin{defn} We say that two m-arc-generators of $L_j$ have the `same type' if both belong to the
same $H_\tau$, or both of them are powers of the same $g_i$ (cf.\ Lemma \ref{lem:z2}).
An m-arc-generator of a certain type is called `pure' if it is not equal
to another m-arc-generator of different type.
\end{defn}
\begin{thm}\label{th:2C} Consider the presentation of $\pi_{1,j}$ from
 (\ref{say:pi1j}) and the m-arc-generators of $L_j$, represented in $\pi_{1,j}$ as in Lemma
 \ref{lem:z2}.
Then if two  m-arc-generators are conjugate in $\pi_{1.j}$ then they are equal in $\pi_{1,j}$.

A non-pure m-arc-generator of  $\pi_{1,j}$  is a power of $h$.
\end{thm}
Proof. We analyze several cases according to the types of the m-arc-generators.
Cases 3 and 4 compare generators of different types.

\vspace{1mm}

{\bf Case 1.} \ Assume that $g_i^r\sim g_i^s$ ($1\leq i\leq n; \ r,s\in\z$).

Then $g_i^{\alpha_ir}\sim g_i^{\alpha_is}$, or $h^r\sim h^s $. Since $h$ is central
and of infinite order, $r=s$.

\vspace{1mm}

{\bf Case 2.} \ Assume that $A\sim A'$, where $A, A'\in H_\tau$.

Let us write $A=h^{r}f_\tau^{c}$ and $A'=h^{s}f_\tau^{d}$.
Since $h$ is central we can assume $s=0$. Set $y:=\tprod_{\sigma>\tau}f_\sigma\cdot
\ell\cdot\tprod_{\sigma<\tau}
f_\sigma$.  Then $h^b=\ell \tprod_{\sigma<\tau}f_\sigma \cdot f_\tau\cdot \tprod_{\sigma>\tau}f_\sigma=
yf_\tau$, hence
$A=y^{-c}h^{cb+r}$ and $A'=y^{-d}h^{db}$.
Projected  into $\overline{\pi_{1,j}}$, we get $\bar{y}^c\sim \bar{y}^d$.
 By the `Conjugacy theorem of free products' \cite[Theorem 1.4, Ch. IV]{LS},  $c=d$.
Next,
write $A'$ as $x^{-1}Ax$. Then $h^r=[x,f^c_\tau]$, hence $h^r$ projected into
$H_1(L_j)=\pi_{1,j}/[\pi_{1,j},\pi_{1,j}]$ is trivial. But the image of $h$ in $H_1(L_j)$
has infinite
order (use the abelianized $G$ from (\ref{say:pi1j}.1)), hence $r=0$ too.

 \vspace{1mm}

{\bf Case 3.} \ Assume that $g_i^r\sim g_j^s$ ($1\leq i,j \leq n, \ i\not= j;  \ r,s\in\z$),
or $g_i^r\sim A$ ($1\leq i\leq n;  \ r\in\z$), where $A\in H_\tau$ ($1\leq \tau\leq k$).

Consider the projection of the two elements in $\overline{\pi_{1,j}}$:
$g_i^r$ projects into the factor  $\z_{\alpha_i}$, while  $H_\tau$ in
$F_{2g+k-1}$. (Hence, either they belong to different factors, or,  if $2g+k=1$, then
$A$ projects into a power of $\ell$ in $*_i\ \z_{\alpha_i}$.)
In any case, the projections  can be conjugate in $\overline{\pi_{1,j}}$ only if they are both trivial,
cf.\  \cite[Theorem 1.4, Ch IV]{LS}.
That means that the two original elements belong to the central subgroup $\langle h\rangle$.
Since $h$ is central, if they are conjugate they should be equal.

\vspace{1mm}

{\bf Case 4.} \ Assume that $A_\tau\sim A_\sigma$, where $A_\tau\in H_\tau$,
$A_\sigma\in H_\sigma$, $\tau\not= \sigma$ ($k\geq 2$).

Let us write $A_\tau=h^{r}f_\tau^{c}$ and $A_\sigma=h^{s}f_\sigma^{d}$. Their  projections
$\bar{f}_\tau^{c}$ and $\bar{f}_\sigma^{d}$  belong
to the factor  $F_{2g+k-1}$
of $\overline{\pi_{1,j}}$.
Assume $2g+k>2$.  When we define the factor
$F_{2g+k-1}$ we can eliminate a variable different from
$f_\tau$ or $f_\sigma$, hence we can write $F_{2g+k-1}$
as  a free product such that
$\bar{f}_\tau^{c}$ and $\bar{f}_\sigma^{d}$ belong to two different non-trivial factors.
Therefore $\bar{f}_\tau^{c}=\bar{f}_\sigma^{d}=1$.
If $F_{2g+k-1}=\z$, we can assume that  its generator is $\bar{f}_\tau$, and
$\bar{f}_\sigma=(\ell\bar{f}_\tau)^{-1}$. Hence $\bar{f}_\tau^{c}\sim
(\ell\bar{f}_\sigma)^{-d}$. Then, again by \cite[Theorem 1.4, Ch IV]{LS}, $c=d=0$.
Hence, in any case, $A_\tau$ and $A_\sigma$  belong
to   $\langle h\rangle$, and they are equal.
 \qed

\begin{say}
Now, using (\ref{lem:marc}) and Theorem \ref{th:2C}
 we are ready to apply Theorem \ref{th:conj} for m-arc and arc-generators.
For $\pi_1$ we  use the notations of (\ref{say:pi1})--(\ref{say:arcgens}).

\end{say}
\begin{thm}\label{th:conj2}
Suppose that $(\sC,\mu)=(v_0,\mu_0)$
 and $(\sC',\mu')=(v_0',\mu_0')$ are two cyclically reduced forms of length zero, whose labels
 $\mu_0$ and $\mu_0'$ are conjugate in $\pi_1$.
 \begin{enumerate}
 \item If $\mu_0$ and $\mu_0'$ are m-arc-generators, then the following cases can happen:
 \begin{enumerate}
 \item (m=0) \ $v_0=v_0'$ and  $\mu_0=\mu_0'$ in $G(v_0)$,
 \item (m=1) \  $v_0$ and $v_0'$ are connected by the edge $e_{\tau_0}$, $\mu_0'=\phi_{e_{\tau_0}}(\mu_0)$,
 and at least one of the generators $\mu_0\in G(e_{\tau_0})^-$
 and $\mu_0'\in G(e_{\tau_0})^+$ is pure.
 \item (m=2) \
 there exists a path
 $(v_0,e_{\tau_0},v_{j_1}, e_{\tau_{1}},v_{0}')$ as in (\ref{th:conj})(2), with $ e_{\tau_{1}}\not=
\overline{ e_{\tau_0}}$, and an element $\mu_1\in G(v_{j_1})$ such that
$\mu_0=\phi_{e_{\tau_0}}^{-1}(\mu_1)$, $\mu_0'=\phi_{e_{\tau_1}}(\mu_1)$,
and both $\mu_0\in G(e_{\tau_0})^-$ and $\mu_0'\in G(e_{\tau_1})^+$ are pure.
 \end{enumerate}
 \item If $\mu_0$ and $\mu_0'$ are arc-generators, then the cases (m=1) and (m=2) from (1) can not happen.
In the (m=0) case the coincidence
$\mu_0=\mu_0'$ of the labels in $\pi_1$ associated with two different arc-generators can happen
only if both are supported by the same chain without arrow (including the generators of type
$h^n$): all of them are  positive powers of the corresponding $g_i$.
 \end{enumerate}
\end{thm}

Proof. (1) We apply Theorem \ref{th:conj}. By Lemma \ref{lem:z2} all elements of type
$\phi^\pm _{\tau_i}(c_i)$  are m-arc-generators in the corresponding
groups $G(v_{j_i})$. Hence, by Theorem \ref{th:2C}, all conjugacies  in Theorem
\ref{th:conj}(1)-(2) are equalities.
Let us write $\mu_i:=\phi^+_{\tau_{i-1}}(c_{i-1})=\phi^-_{\tau_{i}}(c_{i})$.
Then $\phi_{e_{\tau_i}}(\mu_i)=\mu_{i+1}$ for $0\leq i\leq m-1$ (where $\mu_m:=\mu_0'$).
If $e_{\tau_i}=\overline {e_{\tau_{i-1}}}$ is a `backtracking', then $\mu_{i+1}=\mu_{i-1}$,
and the segment of path $(v_{j_{i-1}},e_{\tau_{i-1}},v_{j_i}, e_{\tau_{i}},v_{j_{i+1}})$
can be shortened to $v_{j_{i-1}}$ (although, in such a case,
$\mu_{i+1}$ is special, it is in the image of  $\phi_{e_{\tau_i}}$, and this information
disappears after we shorten the path).

Hence, we can assume that we have no backtracking in the path.

Assume that $m\geq 3$, and we consider the segment of path $
(e_{\tau_0},v_{j_{1}},e_{\tau_{1}},v_{j_2}, e_{\tau_{2}})$.  Then the very same element
$\mu_1:=\phi^+_{\tau_{0}}(c_{0})=\phi^-_{\tau_{1}}(c_{1})$ is represented in two different ways as
m-arc-generator, of  different types (since $\tau_1\not=\overline{\tau_0}$).
 Hence, by Theorem \ref{th:2C}, $\mu_1$ is a power
$h_1^{n_1}$ of the central element $h_1$ of the Seifert piece $G(v_{j_1})$.
The same fact is true in $G(v_{j_2})$, one has  $\mu_2=h_2^{n_2}$.
But, we claim that an identity of type ($\dagger$)
$\phi_{e_{\tau_1}}(h_1^{n_1})=h_2^{n_2}$ cannot  happen unless
$n_1=n_2=0$.  This basically follows from the fact that $G(e_{\tau_1})=\z^2$ injects into $\pi_1$, and
($\dagger$)  would give a relation in this $\z^2$. The proof imitates the proof of Lemma \ref{lem:z2}.
Consider the chain connecting $v_{j_1}$ and $v_{j_2}$, and assume that the Euler numbers are
$-b_1,\ldots, -b_s$. Let $\g_1$ be the arc-generator associated with the vertex adjacent to $v_{j_1}$.
Then $h_2=\phi_{e_{\tau_1}}(\g_1^{\det[b_1,\ldots,b_s]}h_1^{-\det[b_2,\ldots,b_s]})$.
Moreover, $h_1$ and $\g_1$ generate $\z^2$ in $\pi_1$.
But $\det[b_1,\ldots,b_s]\not=0$ (since the graph is negative definite),
hence ($\dagger$) cannot happen.

All the other restrictions in part (1) follow similarly.

(2) Let us start with case {\it (m=1)}, and we assume that $\mu_0\not\in \langle h\rangle$.
 Connect   $v_0$ with $v_0'$ by a chain
with decorations $-b_1,\ldots, -b_s$. Let $h=\g_0, \g_1, \ldots, \g_s, h'=\g_{s+1}$
be the corresponding arc-generators. Fix two integers $0\leq i<j\leq s+1$ and write
$\mu_0=\g_i^{n_i}\g_{i+1}^{n_{i+1}}$, where $n_i\geq 0$ and $n_{i+1}>0$, and
$\mu_0'=\g_j^{m_j}\g_{j+1}^{m_{j+1}}$. The case $\mu_0'\not\in \langle h'\rangle$
 (when both generators are associated with the interior of the chain)
 corresponds to $m_j\geq 0$ and $m_{j+1}>0$, while $\mu_0'\in \langle h'\rangle$
is covered by $m_j=0$, $m_{j+1}\in \z$ for $j=s$.
(The last condition covers the case when $\mu_0'$ is not pure too.)  Then
$$\begin{pmatrix} \ \ [b_{i+1},\ldots, b_j] & \ \ [b_{i+1},\ldots, b_{j-1}]\\
-[b_{i+2},\ldots, b_j] & -[b_{i+2},\ldots, b_{j-1}]\end{pmatrix}\cdot
\begin{pmatrix} m_{j+1}\\ m_j \end{pmatrix}=\begin{pmatrix} n_{i+1}\\ n_i \end{pmatrix}.
$$
But this system has no solution with the above (positivity) restrictions.

The cases  {\it (m=2)}, or {\it (m=0)} of different types, are
 eliminated similarly as the {\it (m=1\,:\, $\mu_0'\in \langle h'\rangle$)} case.
The case {\it (m=0)} of the same $H_\tau$ is eliminated as
{\it  (m=1\,:\,} `interior chain') case.
\qed

\begin{cor} Theorem \ref{sharc.surf.thm.2} is true if $L$ is not a Seifert manifold or a cusp link.
\end{cor}
Proof. This follows from the list of components from Paragraph \ref{sharc.dlt.pairs.say} and
from Theorem \ref{th:conj2}. The coincidences from Theorem \ref{th:conj2}(2)
are eliminated by considering the minimal  dlt modification.  The
powers $g_i^{m_i}$ which cannot be represented as powers of $h$ are exactly those for which
 $m_i/\alpha_i\not\in\z$, cf.\ type (3) in (\ref{sharc.dlt.pairs.say}) (where $m=\alpha_i$). \qed

\begin{say}
 The above proof can be adapted to show Theorem \ref{sharc.surf.thm.2} for cusp singularities as well.
A more precise description is given in   Section \ref{sec.cusps}.
\end{say}

\begin{say}[The Seifert case]\label{say:Seifert} We run the same strategy as above.
The plumbing graph has no arrowheads, and $\pi_1=\pi_1(L)$  and
$\overline{\pi_1}:=\pi_1/\langle h\rangle$  are the following:
\begin{equation*}
\pi_1=\langle h, \{g_i\}_{i=1}^n,\{ a_m,b_m\}_{m=1}^g \
|\ h \ \mbox{central},\ g_i^{\alpha_i}=h,\forall i, \
h^b = \tprod_m [a_m,b_m]\cdot\tprod_i g_i^{\omega_i}\rangle.
\end{equation*}
\begin{equation*}
\overline{\pi_1}=\langle \{\bar{g}_i\}_{i=1}^n,\{ a_m,b_m\}_{m=1}^g \
|\ \bar{g}_i^{\alpha_i}=1,\forall i, \
\tprod_m [a_m,b_m]\cdot\tprod_i \bar{g}_i^{\omega_i}=1\rangle.
\end{equation*}
By  $x_i:=\bar{g}_i^{\omega_i}$, $\bar{g}_i:=x_i^{\omega'_i}$, where
$\omega_i\omega_i'\equiv 1 \mod \alpha_i$, $\overline{\pi_1}$
transforms into the {\it crystallographic group}
(if $g=0$, it is  also named  {\it Dyck--Schwartz polygonal group})
 \begin{equation*}
\langle \{x_i\}_{i=1}^n, \{ a_m,b_m\}_{m=1}^g \
|\ x_1^{\alpha_1}=\cdots =x_n^{\alpha_n}=\tprod_m [a_m,b_m]\cdot \tprod_i x_i=1\rangle.
\end{equation*}
Then $\pi_1$ is infinite iff $\overline{\pi_1}$\, is infinite. Moreover,
infinite crystallographic groups have the following property, see e.g. \cite[page 62]{JN}.
\end{say}
\begin{prop}\label{prop:Sei}
If $\bar{g}_i^r\sim \bar{g}_j^s$ in $\overline{\pi_1}$ ($1\leq i\not=j \leq n;  \, r,s\in\z$),
then $\bar{g}_i^r=\bar{g}_j^s=1$. \qed
\end{prop}
\begin{cor}\label{cor:Sei}
 Theorem \ref{sharc.surf.thm.2} is true if $L$ is a Seifert manifold.
\end{cor}
Proof. Quotient singularities were already treated in Section \ref{sec.quot}.
Otherwise, 
by Proposition \ref{prop:Sei},
 the following arc-generators in $\pi_1$ belong to different
conjugacy classes: $\{h^m\}_{m\in \z_{>0}}$, $\{g_i^{m_i}\}_{m_i\in\z_{>0}\setminus \alpha_i\z}$,
$1\leq i\leq n$. Considering the minimal  dlt modification we hit exactly the
components from Paragraph \ref{sharc.dlt.pairs.say}.
\qed
\begin{say}[Proof of Proposition \ref{quasihomog.prop}]\label{say:uj}
(1) follows from the fact that the classes $h$ and $g_i$ can be realized
by equivariant arcs (parametrizations of the generic and special orbits). (2)  follows
from Corollary \ref{cor:Sei} (or its proof), saying that the conjugacy classes
$\{h^m\}_{m\in \z_{>0}}$, $\{g_i^{m_i}\}_{m_i\in\z_{>0}\setminus \alpha_i\z}$
 are all different. The class $h^m$ can be realized by a family of equivariant arcs
parametrized by the smooth part of the central curve.
\end{say}

\section{Cusp singularities}\label{sec.cusps}

\begin{defn}[Cusps and their  canonical representation]\label{say.cusp.1}
Cusp normal surface singularities are characterized by the fact that the dual graph
of their minimal log resolution is a cyclic graph with all genus decorations zero.
If $E$ is the exceptional curve, then $\pi_1(E)=\z$, and fixing the generator of
$\pi_1(E)$ is equivalent with fixing the cyclic order of the vertices of the graph.
Let us fix such an ordering.
The corresponding Euler numbers of the ordered vertices
will be denoted by $-b_1,\ldots, -b_k$ (up to cyclic permutation), where
each $b_i\geq 2$ and $\sum_i(b_i-2)>0$.

The dual graph $\Gamma$ has no nodes, nevertheless it has exactly one JSJ torus $T$,
which can be chosen as the torus corresponding to one of the edges of the graph.
We rename the vertices in such a way that the cutting edge is $(k,1)$, hence
if we cut the link, or the graph $\Gamma$, along this torus, we get

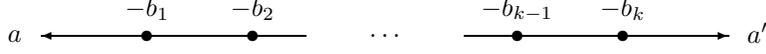
\begin{picture}(300,35)(-50,0)
\put(20,10){\circle*{4}} \put(60,10){\circle*{4}}
\put(160,10){\circle*{4}} \put(200,10){\circle*{4}}
\put(20,10){\line(1,0){60}} \put(140,10){\line(1,0){60}}
\put(20,10){\vector(-1,0){40}} \put(200,10){\vector(1,0){40}}
\put(20,20){\makebox(0,0){$-b_1$}}
\put(60,20){\makebox(0,0){$-b_2$}}
\put(160,20){\makebox(0,0){$-b_{k-1}$}}
\put(200,20){\makebox(0,0){$-b_k$}}
\put(110,10){\makebox(0,0){$\ldots$}}
\put(-30,10){\makebox(0,0){$a$}}
\put(251,11){\makebox(0,0){$a'$}}
\end{picture}

The corresponding arc generators of this string with two arrowheads will be denoted by $
\g_0=\g_a,\g_1,\ldots, \g_k, \g_{k+1}=\g_{a'}$.
By a convenient choice of the connecting paths
 we may assume that they commute with each other and
satisfy the relations
$$\g_i^{b_i}=\g_{i-1}\g_{i+1}, \ \ \ \ \mbox{for $1\leq i\leq k$.}$$
In particular, $\g_0$ and $\g_1$ generate freely $\z^2\simeq\pi_1(T)$, and all
$\g_i$ can be expressed in terms
of them. It is convenient to organize this as follows.

For any collection of integers $a_1,\ldots, a_n$ we define the matrix
$$M(a_1,\ldots,a_n):=\begin{pmatrix} a_1 & 1\\-1& 0\end{pmatrix}\cdots
\begin{pmatrix} a_n & 1\\-1& 0\end{pmatrix}\in {\rm SL}(2,\z).$$
Introduce the column vectors  $\{v_i\}_{i=0}^{k+1}$ of $\z^2$  by $v_0=\binom{0}{1}$ and
$M(b_1,\ldots, b_i)=(v_{i+1},v_i)$ for $1\leq i\leq k$.
E.g., $v_1=\binom{1}{0}$, $v_2=\binom{b_1}{-1}$. For any $v=\binom{r}{s}$ we set
$\underline{\g}^v:=\g_1^r\g_0^s\in \pi_1(T)$. Then, by induction,
$\g_i=\underline{\g}^{v_i}$ for any $0\leq i\leq k+1$.

The matrix $M=M(b_1,\ldots,b_k)\in {\rm SL}(2,\z)$ is called the monodromy operator
associated with the choice of the torus $T$ and the bases $\{v_0,v_1\}$ of $\z^2$,
and the representation  $\g_0\mapsto v_0$, $\g_1\mapsto v_1$ of $\pi_1(T)\to \z^2$.
Note that $Mv_0=v_k$ and $Mv_1=v_{k+1}$. The sequence of
vectors $v_i$ can be extended to a sequence $\{v_i\}_{i\in\z}$ of $\z^2$ by
$v_{\ell k+i}:=M^\ell v_i$ ($\ell\in \z$),
which evidently correspond to arc generators of the   infinite
string $\ldots, b_1,\ldots, b_k, b_1,\ldots, b_k, \ldots$.
Then, the link is a torus bundle over $\s^1$, (where the fiber is identified with $T$), hence
one has the (HNN) extension
$$1\to \z^2=\pi_1(T)\to \pi_1\to  \z=\pi_1(E) \to 1,$$
where $\pi_1(E)$ acts on $\pi_1(T)$
via the monodromy $M$. In other words,
 $$\pi_1=\langle v,t\ |\ v\in\z^2,\ tvt^{-1}=Mv\ \rangle.$$
 This identification of $\pi_1(T)$ with $\z^2$ will be called its `canonical representation'.
 \end{defn}

 \begin{say}[The monodromy operator $M$]\label{say.cusp.2}
 The monodromy operator $M\in{\rm SL}(2,\z)$ satisfies  $\tau :={\rm trace }(M)\geq 3$.
 This shows that the characteristic polynomial $\lambda^2-\tau\lambda+1$
 has two positive (non-rational) roots. We denote them $\lambda_1>1>\lambda_2>0$.
 Let $V_1$ and $V_2$ in $\r^2$ be the corresponding eigenvectors chosen such that $V_1$
 is in the fourth quadrant and $V_2$ in the second one. Then one has the following limits of halflines
 $$\lim_{i\to\infty} \r_{> 0} v_i= \r_{> 0} V_1,\  \
 \lim_{i\to -\infty} \r_{> 0} v_i= \r_{>0} V_2.
 $$
Note that   $( \r_{>0} V_1\cup \r_{>0} V_2)\cap \z^2=\emptyset$.
\end{say}
\begin{prop}\label{prop.cusp.0}
 Let ${\rm Cone}(\Gamma)$ be the real (open) positive cone $\r_{>0}\langle V_1,V_2
\rangle$ generated by $V_1$ and $V_2$, and set
${\rm Cone}_\z(\Gamma):={\rm Cone}(\Gamma)\cap \z^2$.
 Consider also the action of $\z$
on ${\rm Cone}_\z(\Gamma)$ given by $\z\times \z^2\ni (\ell,v)\mapsto M^\ell v$.
Then the image of
$w:\pi_0\bigl(\sharc(0\in X)\bigr)\to
\pi_1\bigl(\link(0\in X)\bigr)/(\mbox{\rm conj})$
is in $\pi_1(T)/(\mbox{\rm conj})$, and under the above `canonical representation'
it is identified with  ${\rm Cone}_\z(\Gamma)/\z={\rm Cone}_\z(\Gamma)/M$.
\end{prop}
Proof. Arc generators of type $\g_i^m$ ($m\in \z_{>0}$) are represented on the lattice points on the
ray $\r_{>0}v_i$, and those of type $\g_i^{m_i}\g_{i+1}^{m_{i+1}}$ ($m_i,m_{i+1}\in\z_{>0}$)
by the lattice points
in the open cone    $\r_{>0}\langle v_i,v_{i+1}\rangle$, modulo the monodromy action.
\qed
\begin{say}[Dependencies]\label{say.cusp.3}
$M$  depends on the above choice of `canonical
 representation' (of  $\pi_1\to \z^2$, $\g_0\mapsto v_0$, $\g_1\mapsto v_1$).
If  the JSJ torus is given by the edge $(1,2)$, then the new monodromy
operator is $M_{(1,2)}=M(b_2,\ldots, b_k,b_1)=M(b_1)^{-1}\cdot M\cdot M(b_1)$, hence
$M_{(1,2)}$ and $M$ are conjugate. The same is true for any other choice of $T$.

In fact, the correspondence $(b_1,\ldots, b_k)\mapsto M$
  from {\it oriented cycles}  to the set of conjugacy classes of $A\in {\rm SL}(2,\z)$ with
   ${\rm trace}(A)\geq 3$ is a bijection (cf.\ \cite[6.3]{NP}).
(One recovers $(b_1,\ldots, b_k)$ from
$M$ via an infinite periodic continued fraction associated with 
${\rm Cone}_\z(\Gamma)$.)

If we change the orientation in the cycle of the graph (or, the generator of $\pi_1(T)$), then
the new monodromy operator will be $\overline{M}=M(b_k,\ldots, b_1)$, and
$$\overline{M}=SM^{-1}S^{-1}, \ \ \mbox{where } \ \ S=S^{-1}=\begin{pmatrix} 0&1\\1&0\end{pmatrix}.$$
If $T(A)$ denotes the torus bundle over $\s^1$ with monodromy $A$ (where ${\rm trace}(A)\geq 3$
always), then $T(A)$ is orientation preserving diffeomorphic with $T(B)$ iff $A$ is conjugate in
${\rm SL}(2,\z)$ to either $B$ or $SB^{-1}S^{-1}$ \cite[6.2]{NP}.
Hence $T(M)\simeq T(\overline{M})$, and their fundamental groups are isomorphic via
 $t\leftrightarrow t^{-1}$, $v\leftrightarrow Sv$  and $M\leftrightarrow M^{-1}$.

Changing the orientation of $\partial \bdd$, the image of $w$ is the set of
lattice points of  $-{\rm Cone}_\z(\Gamma)$ up to the monodromy action of $M$.
Our next goal is to identify the lattice points of the complementary cone
of ${\rm Cone}(\Gamma)$ (up to the action of $M$)
with   $\pi_0\bigl(\sharc(0\in X^*)\bigr)$, where $X^*$ is the dual cusp.
\end{say}
 \begin{say}[The dual cusp]\label{say.cusp.4}
 We wish to identify $\pi_1(T)$ for the cusp $X$ and its dual $X^*$ in a canonical way, such that
 both monodromy actions will act on the same $\z^2$.

 A possible plumbing graph of the link $L^*$ of the dual cusp  $X^*$ is $-\Gamma$. This means that we replace
 each decoration $-b_i$ of $\Gamma$  by $b_i$, and each edge decoration (which was  $+$) by $-$.
 The effect of this is that $L^*$  is $-L$, $L$ with opposite orientation.

  By (oriented) plumbing calculus (see \cite{NP}) $-\Gamma$ can be replaced by another graph, $\Gamma^*$,
  which has all Euler decorations $\leq -2$ and edge decorations $+$, and it is the dual graph
  of the minimal log resolution of $X^*$. Before we describe it,
  we fix/identify the cutting JSJ tori in both 3-manifolds.
  Let us fix a cyclic ordering  and the cutting edge in $\Gamma$ as in (\ref{say.cusp.1});
  and mark the corresponding Euler numbers by $-b_1,\ldots, -b_k$.
  Since at least one of them is $\leq - 3$, we can assume that $b_k\geq 3$.

  We choose in $-\Gamma$ the very same edge as cutting edge and the same orientation of the cycle
  (generator of $\pi_1(E)$). Then, after running the plumbing calculus (in such a way that we never
  blow up/down the cutting edge) we get the vertices with decorations $-b_1^*,\ldots , -b^*_{k^*}$.
 One checks (see e.g. \cite{NP}) that if the $b$-sequence is
 $$2^{k_1^*-1}, k_1+2, 2^{k_2^*-1}, k_2+2,\ldots, k_g+2,$$
 then the $b^*$ sequence in $\Gamma^*$ is
    $$k_1^*+2, 2^{k_1-1}, k_2^*+2, 2^{k_2-1}, \ldots, 2^{k_g-1},$$
    where $2^k$ represents $k$ copies of 2.

 The goal is to identify the two tori-fibers, that is, to see both tori in the same graph. This seems to appear only implicitly in the literature
 \cite{NP, PPP}.
 Consider the following plumbing graph:

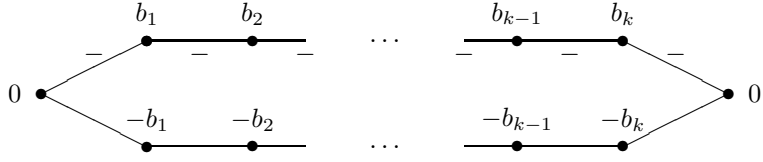
\begin{picture}(300,75)(-50,0)
\put(20,50){\circle*{4}} \put(60,50){\circle*{4}}
\put(160,50){\circle*{4}} \put(200,50){\circle*{4}}
\put(20,50){\line(1,0){60}} \put(140,50){\line(1,0){60}}
\put(20,50){\line(-2,-1){40}} \put(200,50){\line(2,-1){40}}
\put(20,60){\makebox(0,0){$b_1$}}
\put(60,60){\makebox(0,0){$b_2$}}
\put(160,60){\makebox(0,0){$b_{k-1}$}}
\put(200,60){\makebox(0,0){$b_k$}}
\put(110,50){\makebox(0,0){$\ldots$}}
\put(-30,30){\makebox(0,0){$0$}}
\put(250,30){\makebox(0,0){$0$}}

\put(40,45){\makebox(0,0){$-$}}\put(0,45){\makebox(0,0){$-$}}
\put(80,45){\makebox(0,0){$-$}}
\put(140,45){\makebox(0,0){$-$}}
\put(180,45){\makebox(0,0){$-$}}
\put(220,45){\makebox(0,0){$-$}}

\put(20,50){\circle*{4}} \put(60,10){\circle*{4}}
\put(160,10){\circle*{4}} \put(200,10){\circle*{4}}
\put(20,10){\line(1,0){60}} \put(140,10){\line(1,0){60}}
\put(20,10){\line(-2,1){40}} \put(200,10){\line(2,1){40}}
\put(20,20){\makebox(0,0){$-b_1$}}
\put(60,20){\makebox(0,0){$-b_2$}}
\put(160,20){\makebox(0,0){$-b_{k-1}$}}
\put(200,20){\makebox(0,0){$-b_k$}}
\put(110,10){\makebox(0,0){$\ldots$}}

\put(20,10){\circle*{4}}\put(-20,30){\circle*{4}}\put(240,30){\circle*{4}}
\end{picture}

By plumbing calculus (0-chain and oriented handle absorptions), this graph is equivalent
to  a single vertex with genus 1 and Euler number 0, hence representing
 the trivial torus bundle over $\s^1$. Since the monodromy acting on the torus is trivial, the tori associated with the edges can be identified by canonical
isomorphisms.

In this graph we can proceed on the `top subchain' the calculus which provided the $b^*$ sequence of $\Gamma^*$
(just repeating the steps which provides the normal form of $-\Gamma$). Then we get the graph
(here we use the fact that $b_k\geq 3$):

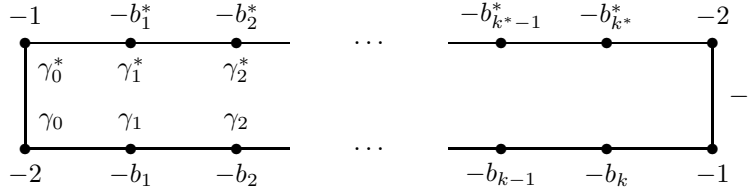
\begin{picture}(300,80)(-50,-5)
\put(20,50){\circle*{4}} \put(60,50){\circle*{4}}
\put(-20,50){\circle*{4}} \put(240,50){\circle*{4}}
\put(-20,10){\circle*{4}} \put(240,10){\circle*{4}}
\put(160,50){\circle*{4}} \put(200,50){\circle*{4}}
\put(20,50){\line(1,0){60}} \put(140,50){\line(1,0){60}}
\put(-20,50){\line(0,-1){40}} \put(240,50){\line(0,-1){40}}

\put(20,50){\line(-1,0){40}} \put(200,50){\line(1,0){40}}

\put(-10,40){\makebox(0,0){$\g^*_0$}}
\put(20,40){\makebox(0,0){$\g^*_1$}}
\put(60,40){\makebox(0,0){$\g^*_2$}}
\put(-10,20){\makebox(0,0){$\g_0$}}
\put(20,20){\makebox(0,0){$\g_1$}}
\put(60,20){\makebox(0,0){$\g_2$}}

\put(20,60){\makebox(0,0){$-b^*_1$}}
\put(60,60){\makebox(0,0){$-b^*_2$}}
\put(160,60){\makebox(0,0){$-b^*_{k^*-1}$}}
\put(200,60){\makebox(0,0){$-b^*_{k^*}$}}
\put(110,50){\makebox(0,0){$\ldots$}}
\put(250,30){\makebox(0,0){$-$}}

\put(-20,0){\makebox(0,0){$-2$}}
\put(-20,60){\makebox(0,0){$-1$}}
\put(240,60){\makebox(0,0){$-2$}}
\put(240,0){\makebox(0,0){$-1$}}

\put(20,10){\circle*{4}} \put(60,10){\circle*{4}}
\put(160,10){\circle*{4}} \put(200,10){\circle*{4}}
\put(20,10){\line(1,0){60}} \put(140,10){\line(1,0){60}}
\put(20,10){\line(-2,0){40}} \put(200,10){\line(2,0){40}}
\put(20,0){\makebox(0,0){$-b_1$}}
\put(60,0){\makebox(0,0){$-b_2$}}
\put(160,0){\makebox(0,0){$-b_{k-1}$}}
\put(200,0){\makebox(0,0){$-b_k$}}
\put(110,10){\makebox(0,0){$\ldots$}}
\end{picture}

Only one edge has decoration $-$, the others have $+$. We inserted also some of the
arc generators associated with the vertices compatibly to the notation of  (\ref{say.cusp.1}).

Next, we fix $\pi_1(T)=\z^2$ under the canonical  identification of $\Gamma$ as in  subsection (\ref{say.cusp.1})
with base elements $v_0=\binom{0}{1}$ and $v_1=\binom{1}{0}$, and we represent $\gamma_i=\underline{\g}^{v_i}$
for $i=1,2$. Then, we  express in the above graph the elements $\g^*_i$ ($i=1,2$) in terms of $\g_0$ and $\g_1$.
Namely, using the relations $\g_0^2=\g_1\g^*_0$ and $\g_0^*=\g_0\g^*_1$, we get
$\g^*_0=\g_0^2\g_1^{-1}$ and $\g_0\g_1^{-1}$. Hence, although in the canonical representation of $\Gamma^*$ (considered independently from $\Gamma$, and repeating (\ref{say.cusp.1}) for it),
we would send  $\g^*_0$ to $\binom{0}{1}$ and $\g^*_1$ to $\binom{1}{0}$,
in this representation compatible with the canonical representation of $\Gamma$
(when  we send  $\g_0$ to $\binom{0}{1}$ and $\g_1$ to $\binom{1}{0}$)
we have to send $\g^*_0$ to $u_0:=\binom{-1}{2}$ and $\g^*_1$ to $u_1:=\binom{-1}{1}$.

Let ${\rm Cone}^c(\Gamma)$ be the complementary cone of ${\rm Cone}(\Gamma)$, namely
$\r_{>0}\langle -V_1, V_2\rangle$, and set also ${\rm Cone}^c_\z(\Gamma)={\rm Cone}^c(\Gamma)\cap \z^2$.
Then using the properties of $M$ one verifies that $u_0,u_1\in {\rm Cone}^c_\z(\Gamma)$.
 Moreover, one has the following identity;
cf.\ \cite[7.6]{nakamura}.  \end{say}

 \begin{prop}\label{prop.cusp.1}
 Set $T:=\begin{pmatrix} -1&-1\\ 1&2\end{pmatrix}$
 and consider the monodromy operators
 $M=M(b_1,\ldots, b_k)$ and $M^*=M(b_1^*,\ldots, b^*_{k^*})$. Then $MT=TM^*$,
hence the eigenvalues of $M$ and $M^*$  are the same.
 \end{prop}
 Proof. In the above graph the monodromy operator of the trivial torus-fibration is the identity.
 This reads  as (starting from the left $-1$-vertex)
 $$M(1,2)\cdot M\cdot M(1,2)\cdot M(b^*_{k^*}, \ldots, b^*_1)=-I,$$
 where the sign is given by the unique negative edge. Since $M(1,2)=ST^{-1}=-TS$ and by (\ref{say.cusp.3})
 one also has $M(b^*_{k^*}, \ldots, b^*_1)=S(M^*)^{-1}S$, we get $MT=TM^*$.
 \qed

 \vspace{2mm}

 Note that $\det(T)=-1$, hence  the  proposition does not say that $M$ and $M^*$ are necessarily conjugate in
 ${\rm SL}(2,\z)$ (although, in special cases they can be conjugate or even equal, see e.g the case
 $(3,3,3)$ appearing in Example \ref{cusp.exmp.i}, which is auto-dual).

 \begin{cor}\label{cor.cusp.1} 
Let $V_1^*$ and $V_2^*$ be the eigenvectors of $M^*$ defined via the canonical
 recipe of (\ref{say.cusp.2}). Then
\begin{enumerate}
\item   $TV_1^*\in -\r_{>0}V_1$ and  $TV_2^*\in \r_{>0} V_2$.
 In particular, $T$ sends ${\rm Cone}_\z(\Gamma^*)$ isomorphically onto
 ${\rm Cone}_\z^c(\Gamma)$ identifying the action of 
$M^*$ on the first one with the action of
 $M$ on the second  one. 
\item  ${\rm Cone}_\z^c(\Gamma)/M $ is identified via $T^{-1}$
 with ${\rm Cone}_\z(\Gamma^*)/M^*$ and Proposition \ref{prop.cusp.0}
identifies the latter with
$\pi_0\bigl(\sharc(0\in X^*)\bigr)$. \qed
\end{enumerate}
 \end{cor}

\begin{rem}[The four cones and orientations]
Let us fix a generator of $\pi_1(E)$ as above. 
There are two sources of orientations on the torus $T$.

First, regarding $L$ as a torus
bundle over $S^1$,
$\pi_1(S^1)$ and $\pi_1(E)$ are canonically identified, hence an
orientation of $L$ is equivalent to an orientation  of $T$.

 The plumbing representation of $L$ gives another orientation as  we
glue several
$S^1$-bundles over  annuli.
Let us call the circle in the  annulus  $[0,1] \times S^1$ the `angle
circle' and the
$S^1$-fiber of the bundles the `fiber circle'.  Then $T$ consists of their
product, and
the orientation of $T$ is the product of angle  and  fiber circle
orientations.
Switching both circle orientations does not alter  orientation of the torus 
or the
link.
The realization of the link as a tubular
neighborhood of holomorphic curves fixes
an orientation of the fiber circle. The latter is not a topological invariant
of the  oriented link. 
\end{rem}

\section{Inoue surfaces and short arcs on cusps }\label{sec.inoue}

Short arcs on a cusp and on its dual can be seen together
very nicely on hyperbolic Inoue surfaces \cite{MR0442296, MR0442297}.
We need to understand mostly their cusps, treated in detail
in \cite[\S.2]{MR0393045}.

\begin{defn}[Hyperbolic Inoue surfaces]
Let $K=\q\bigl(\sqrt{d}\bigr)\subset \r$ be a real quadratic field.
Conjugation is denoted by $w\mapsto w'$.

Let $H\subset K$ be a free $\z$-module of rank 2. Let
$u\in K$ be a nontrivial element such that $u>0, u'>0$ and $uH=H$.
Let $V=\langle u\rangle$ be the group it generates.
We may assume that $H$ is an ideal in $\o_K$; in this case $u$ is a 
unit in $\o_K$ and this makes
some formulas simpler.

Let $G=G(K,H, u)\subset \GL(2,\r)$ be the group 
$$
G:=\Bigl\{
\Bigl(
\begin{array}{cc}
v & m\\
0 & 1
\end{array}
\Bigr) :
v\in V, m\in H\Bigr\}.
$$
There is a natural group extension
$$
0\to H\to G\to V\to 0.
$$
The group $G$ acts  on $\c^2$ by the rule
$$
\rho(v, m): (z_1, z_2)\mapsto  (vz_1+m, v'z_2+m').
$$
The action is properly discontinuous on $\h\times \c$
where $\h\subset \c$ is the upper half plane.
It is proved in \cite[\S.2]{MR0393045} that the quotient
$(\h\times \c)/G$ can be compactified by adding 2
points; resulting in a singular complex surface with 2 cusps.
The sets  
$$
\{(z_1, z_2): \Im z_1\cdot \Im z_2>1/\epsilon\}\qtq{resp.} \
\{(z_1, z_2): \Im z_1\cdot \Im z_2<-1/\epsilon\}
$$ 
give open neighborhoods of
the $+\infty$-cusp (resp.\ $-\infty$-cusp).

The  minimal resolution of this singular surface is  called a
{\it hyperbolic Inoue surface.}
\end{defn}

\begin{say}[Topology of the hyperbolic Inoue surfaces]\label{top.of.inoue}
To understand the  above construction
 topologically, consider the imaginary parts map
$$
\Im: \h\times \c\to \r^{>0}\times \r\qtq{given by}
(z_1, z_2)\mapsto (y_1=\Im z_1, y_2=\Im z_2).
$$
This is $G$-equivariant. The subgroup $H$ acts trivially
on $\r^{>0}\times \r $ and $V$ acts by
$$
\rho(v): (y_1, y_2)\mapsto  (vy_1, v^{-1}y_2).
$$
Thus $(\h\times \c)/H$ is a torus bundle over
$\r^{>0}\times \r $. It is convenient to extend the 
$G$-action to $\bar\h\times \c$. The $H$-action is still properly
discontinuous and $(\bar\h\times \c)/H$ is a torus bundle over
$\r^{\geq 0}\times \r $. The advantage is that the fiber
over $(0,0)\in \r^{\geq 0}\times \r $ is canonically identified
with 
$\r\otimes_\q K/H$, giving a    canonical isomorphism
$$
H_1\bigl( \Im^{-1}(0,0)/H, \z\bigr)\cong H.
\eqno{(\ref{top.of.inoue}.1)}
$$
(For other fibers, this identification is possible only up-to
the $V$-action on $H$.)

We see that the real hypersurface
$$
T:=(\Im z_2=0)\subset \h\times \c
$$
is $\rho$-invariant. The quotient  $T/G$ is a torus bundle over $\s^1$.

The complement of $T$ is decomposed into 2 pieces
$$
W^+:=\h\times \h \qtq{and} W^-:=\h\times (-\h).
$$
The $G$-equivariat map $(z_1, z_2)\mapsto (z_1, \Re z_2) $
shows that
$$
W^+/G\sim (T/G)\times \r^{>0}\qtq{and}W^-/G\sim (T/G)\times \r^{<0}.
$$
Thus $W^+/G\cup\{\infty\}$ and $W^-/G\cup\{-\infty\}$
are both homeomorphic to cones over $T/G$.
\end{say}

\begin{say}[Short arcs on hyperbolic Inoue surfaces]\label{sharc.on.inoue}
One can construct representatives of the spaces of short
arcs of the two cusps as follows.

First let $m\in H$ be an element such that $m>0, m'>0$.
Then
$$
\tilde \phi_m:\h\to W^+ \qtq{given by} w\mapsto (mw, m'w)
$$
descends to an holomorphic map
$\phi^*_m: \dd^*\cong \h/\z[1]\into W^+/G$. Here
$\z[1]$ denotes the translation action $w\mapsto w+1$
and the quotient $\h/\z[1]$ is identified with the punctured
unit disc $\dd^*$. 
Next $\phi^*_m $ extends to a short arc
$$
\phi_m: \dd\into W^+/G\cup\{\infty\}.
$$
Similarly, if $m>0$ but $ m'<0$ then we get 
a short arc
$$
\phi_m: \dd\into W^-/G\cup\{\infty\}.
$$
Note that  $\tilde \phi_m $ extends to 
$$
\bar \phi_m:\bar \h\to \bar \h\times \bar \h 
$$
and the image of $\r=\partial \bar \h$ gives a
homology class in $\Im^{-1}(0,0)/H$.
Under the identification (\ref{top.of.inoue}.1), 
this homology class is exactly $m\in  H$.

More generally, any short arc through the $+\infty$ cusp
lifts to  a holomorphic map
$\phi: \h\to \h\times \h$ that satisfies the condition
$$
\phi(w+1)=\phi(w)+(m,m')\qtq{for some} m\in H.
$$
  Then $\psi(w):=\phi(w)-(mw, m'w)$ satisfies
$\psi(w+1)=\psi(w)$.
Thus $\psi$ can be expanded into a Fourier series
$$
\psi(w)=\bigl(\tsum_{n\in \z}a_n e^{2\pi inw}, 
\tsum_{n\in \z}b_n e^{2\pi inw}\bigr).
$$
We need to extend this across the cusp, hence
$\psi$, as a function on $\dd^*=\h/\z[1]$,
 can not have essential singularities. This gives the
general solution
$$
 \phi(w)=\bigl(mw+\tsum_{n\geq 0}a_n e^{2\pi inw}, 
m'w+\tsum_{n\geq 0}b_n e^{2\pi inw}\bigr).
\eqno{(\ref{sharc.on.inoue}.1)}
$$

Note finaly that the monodromy action of Definition \ref{say.cusp.1}
 corresponds to
multiplication by  $u$ on $H$. The 4 cones in Paragraph \ref{say.cusp.2}
are given by the conditions
$$
(m>0, m'>0), \ (m>0, m'<0), \ (m<0, m'>0), \ (m<0, m'<0).
$$
We see that short arcs on the $+\infty$ cusp give
elements in the first cone, and short arcs on the $-\infty$ cusp give
elements in the second cone.
Thus our computations give the same results as the topological
considerations of Section \ref{sec.cusps}.

\end{say}

\begin{rem} It is interesting that writing families of short arcs as in
(\ref{sharc.on.inoue}.1) treats all of them the same way.
By contrast, our original approach  suggests
that arc families arising from {\em irreducible components} of the 
exceptional set of the
minimal resolution and arc families arising from the {\em singularities} of the 
exceptional set do behave differently.

By Paragraph \ref{sharc.dlt.pairs.say},
these two types of families can be distinguished topologically
as long as the irreducible component has either genus $\geq 1$
or has at least 3 nodal points. For cusps, every irreducible component
of the exceptional set is a rational curve with 2 nodal points.
\end{rem}

\section{Open problems}\label{sec.conj}

There are many open questions about arc spaces. 
Here we collect some of them that seem most interesting to us.

\subsection*{Local structure of arc spaces}{\ }

\begin{defn}  Let $X$ be an analytic space.
A {\it finite type} family of arcs on $X$ is given by
\begin{enumerate}
\item an  analytic space  $V\subset \c^r_{\mathbf x}$ and 
\item  a holomorphic vector function  
$$
F\bigl(x_1,\dots, x_r; y_1, \dots, y_s, t\bigr):
V\times \dd^s_{\mathbf y}\times \bdd_t\to X
$$
\end{enumerate}
 such that $t\mapsto F\bigl(c_1,\dots, c_r; u_1(t),\dots, u_s(t),t\bigr) $
is an arc in $X$ for all 
$(c_1,\dots, c_r)\in V$ and 
$u_j(t)\in  B_{<1}\hol(\o_{\bdd})$ for $j=1,\dots, s$.

Thus  $F$ gives a family of arcs parametrized by
$V\times B_{<1}\hol(\o_{\bdd})^s$. 

Note that a finite type family of arcs is much more
restrictive than a family of arcs parametrized by
$V\times B_{<1}\hol(\o_{\bdd})^s$.
\end{defn}

We have not checked the details but the methods of
\cite{MR1320605, gri-kaz, drinf} should imply that
if $\phi:\bdd\to X$ is an arc whose image is not contained in
$\sing X$ then $[\phi]\in \arc(X)$ has an open neighborhood
$[\phi]\in U(\phi)\subset  \arc(X)$ 
that is obtained from a finite type family of arcs.

\begin{conj}\label{local.atlas.conj} Let $X$ be a complex space
and $\arc^*(X)\subset \arc(X)$ the space of arcs that are
not contained in $\sing X$.

Then $\arc^*(X)$ has a natural atlas consisting of
finite type families.
\end{conj}

\begin{rem} We are intentionally vague about what
types of  coordinate  changes
we allow in such a  ``natural atlas.'' 
Our hope is that the following three types of operations suffice.
\begin{enumerate}
\item (Shrinking $V$) Replacing $V$ with an open subset of it.
\item  (Changing coordinates) Biholomorphisms 
$V_1  \cong   V_2$ and
$V_1\times \dd^s  \cong   V_2\times \dd^s $
(the latter defined in neighborhoods of  $V_i\times \{0\}$)
that are compatible with the projections
$$
\begin{array}{ccc}
V_1\times \dd^s & \cong  & V_2\times \dd^s\\
\downarrow && \downarrow\\
V_1 & \cong  & V_2
\end{array}
$$
\item (Splitting off the leading coefficient)
A finite type family $(F:V\times \dd^s\times \bdd_t\to X)$ is replaced by
$\bigl(F^s:(V\times \dd)\times \dd^s\times \bdd\to X\bigr)$
using the rule
$$
\bigl(c_1,\dots, c_r; u_1(t),\dots, u_s(t)\bigr)\mapsto
 \bigl(c_1,\dots, c_r, u_1(0); 
\tfrac{u_1(t)-u_1(0)}{t},u_2(t),\dots, u_s(t)\bigr).
$$
\end{enumerate}
\end{rem}

The methods and results  of \cite{MR1320605, MR2663645}
seem quite close to a solution along the above lines.

Even weaker versions of this would be quite important.
For instance, without a complex structure on $\arc^*(X)$,
it is not even clear how to define the notion of irreducible components.
Our Definition \ref{irred.defn} seems to work, but it leaves
several basic questions unanswered. The most important is  the following.

\subsection*{Curve selection conjecture}{\ }

Let $X$ be an irreducible algebraic variety and
$p_1,\dots, p_m\in X$ a collection of points. Then
there is always an irreducible curve $C\subset X$ that contains
all these points.
This implies that there is an arc $\phi:\dd\to X$ whose image contains
the points $p_1,\dots, p_m$. For irreducible complex spaces
much stronger results are proved in \cite{MR2126216, MR2150882}.

% Presumably the latter still holds for
% irreducible complex spaces, though we have not been able to locate a 
% reference.

A similar assertion is much harder for infinite dimensional
complex spaces, in particular for the arc spaces $\arco(X)$ and
$\sharc(X)$.

\begin{conj}\label{irred.prob} Let $X$ be a complex space, $W$ an irreducible
component of $\arco(X)$ or $\sharc(X)$ and
$[\phi_1], \dots, [\phi_m]$ arcs in $W$. 
Then there is a holomorphic family of arcs
$$
F: \dd\times \bdd \to X\qtq{and} p_1,\dots, p_m\in \dd
$$
such that $\phi_i(t)=F(p_i, t )$ for $i=1,\dots, m$. 
\end{conj}

The notion of a strongly irreducible space was introduced in
Definition \ref{irred.defn} to
go around this conjecture.
We do not know how to prove that a connected component of
$\arc(X)$ is strongly irreducible iff it contains a 
strongly irreducible dense open subset. These types of questions
are quite difficult; see for instance  \cite{2012arXiv1201.6310F}.

\subsection*{Is $\arco(0\in X)$ well defined?}{\ }

One usually thinks of a singularity as an equivalence class
of pointed singular spaces $(0\in X)$. The set of short arcs
$\sharc(X)$ does depend on the choice of the representative
$(0\in X)$ but we saw in Paragraph \ref{conn.comp.say} that
the resulting arc spaces  $\sharc(X)$ are naturally homeomorphic to
each other.

If several points of an arc pass through the singularities then the
rescaling trick used in Paragraph \ref{conn.comp.say} does not work
and  the homeomorphism type of $\arco(X)$ 
does depend on the choice of the representative $(0\in X)$.
For instance,  
$$\arco(0\in \dd)\qtq{and} 
\arco\bigl(0\in \dd\setminus\{\tfrac12\}\bigr)
$$
have different connected components.

A positive answer to the following would be a
suitable replacement.

\begin{conj}\label{indep.of.U.conj}
 Let $(0\in X)$ be an isolated singularity and
$0\in U\subset X$ a contractible, Stein neighborhood. Then
the homeomorphism types of  $\arc(U)$ and of $\arco(U)$
are independent of $U$.
\end{conj}

For singularities with a good $\c^*$-action, sufficiently convex
neighborhoods $U$ do give homeomorphic arc spaces,
but we have very little other evidence to support the above conjecture.

\subsection*{Computing $\arco(X)$ for surface singularities}{\ }

Let $(0\in S)$ be a contractible, Stein, normal surface singularity.
Assuming that at least $\pi_0\bigl(\arco(S)\bigr)$
is independent of  the choice of the representative $S$, 
the natural problem is the following variant of
Theorem \ref{sharc.surf.thm.2}.

\begin{ques}\label{arco.surf.ques}
 Let $(0\in S)$ be a contractible, Stein, normal surface singularity.
Is the winding number map
(\ref{conn.comp.say}.1)
$$
\pi_0\bigl(\arco(S)\bigr)\to \pi_1\bigl(\link(0\in S)\bigr)/
(\mbox{\rm conj})
$$
 an injection? 
\end{ques}

We have not been able to compute any examples.
It is also possible that  a connected component
of $\arco(0\in S) $ has several
irreducible components and determining them may well be the
more important question.

\subsection*{Arcs on complex manifolds}{\ }

It would also be interesting to study arcs on
complex manifolds, even in the normal crossing cases.

Let $X$ be a smooth, affine variety
(or a Stein manifold) and $Z\subset X$ a normal crossing divisor.
As before, let $\arco(Z\subset X)$ denote 
 the space of those arcs $\phi:\bdd\to X$
for which  $\phi(\partial \bdd)\subset X\setminus Z$.
We get  a winding number map
$$
w_{X,Z}:\pi_0\bigl(\arco(Z\subset X)\bigr)\to
\ker\bigl[\pi_1(X\setminus Z)\to \pi_1(X)]/(\mbox{conj}).
$$
Note that the kernel on the right hand side may not even be finitely generated.

There are some cases when $Z\neq 0$ yet $\pi_1(X\setminus Z)=1$.
This happens for instance for the affine quadric
$$
X=(x^2+y^2=z^2+1)\subset \c^3\qtq{and} Z=(x-z=y-1=0).
$$
There are, however, many cases when every irreducible component of
$Z$ adds a new generator to  $\pi_1(X\setminus Z)$. 
This happens if $X=\c^n$ or, more generally, if $H^2(X, \z)=0$.

\begin{prob} Study the above winding number map, or its
restriction to short arcs,  in some interesting 
global situations.
\end{prob}

A starting case could be when $X=\c^n$ and $Z$ is a union of hyperplanes
intersecting transversally.

\subsection*{Short arcs in higher dimensions}{\ }

As the examples in \cite{MR2030097, df-arc, k-nash2} show,
the original Nash problem needs to be reformulated in higher dimensions.
A variant was proposed in \cite{k-nash2}  and 
it is not hard to develop a version 
that applies to $\sharc(0\in X)$ in all dimensions.

However, we feel that both the questions in  \cite{k-nash2}
and their $\sharc$ versions need to be tested in some
nontrivial cases. For $cA$-type singularities, the
irreducible components of $\sharc(0\in X)$ are studied in
\cite{john-kol3}.

\subsection*{Real arcs}{\ }

\begin{defn}[Real arcs]\label{real.arcs.defn}
Let $X$ be  a real analytic space.
A  {\it  real analytic arc} in $X$ is a real analytic  morphism
$\phi:[-1,1] \to X$. 
(As before, $\phi$ is defined and analytic in some neighborhood of $[-1,1]$.)
If, in addition, $\supp\phi^{-1}(\sing X)=\{0\}$ then $\phi$ is called
a {\it short  real analytic arc} or 
{\it short arc.} 

For real analytic arcs one should use the topology defined by 
convergence in all derivatives. 
(In the  real case it makes sense to talk about
 $C^m$ arcs where $m\in \{0,1,\dots,\infty\}$.
 These behave quite differently
from real analytic arcs, even for $m=\infty$.) 
\end{defn}

The basic observation is the following finiteness result.

\begin{prop} \label{sharcr.finite.prop}
The space of short real arcs
 $\sharcr(X)$
has only finitely many irreducible components.
\end{prop}

Proof. Let $f:Y\to X$ be a log resolution with exceptional divisor
$E$.  As in Section \ref{sec.dlt}, it is enough to prove that
$\sharcr(E\subset Y)$ has only finitely many irreducible components.

We follow the arguments in Paragraphs
 \ref{snc.short.arcs.say}--\ref{snc.short.arcs.lem}. 
In the complex case, the irreducible components of $\sharc(E\subset Y)$
were described by the strata of $E$ and by the intersection numbers
with the divisors $E_i\subset E$.

In the real case only the parity of the  intersection numbers  matters
as shown by the deformations
$$
(t,s)\mapsto  t^{2r+1}+s^2t \qtq{and} (t,s)\mapsto  t^{2r}+s^2t^2.
$$
We also need some new information. First, by strata we mean
not the irreducible components of the intersections
but their connected components.
Furthermore, if $E_j$ is 2-sided in $X$ then 
an extra datum that appears is  the side that contains the
image of $(0,1]$.

All together we get  finitely many irreducible components
for $\sharcr(E\subset Y)$ and so finitely many irreducible components
for $\sharcr(X)$. \qed
\medskip

We know very little about the connected or irreducible components
for $\sharcr(X)$ in higher dimensions, but the
surface case should be easier to handle.

\begin{exmp}[$A$-type singularities]\label{real.A.exmp}
 There are two real forms of 
$A$-type singularities with enough real points.

\medskip

\ref{real.A.exmp}.1. $(x^2+y^2=z^m)$ These are quotient singularities. 
If $m$ is even,
one of the exceptional curves of the minimal resolutions is real,
if  $m$ is odd then none of them are real. This shows the following.

 $\sharcr\bigl(x^2+y^2=2z^{2m+1}\bigr)$ has
one connected component. A typical general arc is
$\bigl(t^{2m+1}, t^{2m+1},t^2\bigr)$.

 $\sharcr\bigl(x^2+y^2=2z^{2m}\bigr)$ has
two connected components.  Typical general arcs are
$\bigl(t^{m}, t^{m},\pm t\bigr)$.
\medskip

\ref{real.A.exmp}.2. $(xy=z^m)$ These are not quotient singularities.
All $m-1$  exceptional curves of the minimal resolution are real.
Following the arguments of Proposition \ref{sharcr.finite.prop},
we expect 2 types of  connected components in 
$\sharcr(xy=z^m)$.
\medskip

\ref{real.A.exmp}.2.a.  Typical arcs are  $(t^i, t^{m-i}, t)$
for $0<i<m$. These are the arcs whose lift to the  minimal resolution 
intersects one exceptional curve with multiplicity 1.
Adding signs to two of these terms we get
 $4m-8$ connected components of the space of short real arcs.
(The curves at the end give only 2 components each.)
\medskip

\ref{real.A.exmp}.2.b.  Typical arcs are  $(t^i, t^{2m-i}, t^2)$
for $0<i<2m$. These are the arcs whose lift to the  minimal resolution 
intersects either one exceptional curve with multiplicity 2 (for $i$ even)
or two exceptional curves (for $i$ odd). 

Most of the time, these do not give new components. 
For instance, if $i<m$ then we have the deformation
$$
\bigl(t^i, t^{m-i}(t-\epsilon)^m, t(t-\epsilon)\bigr),
$$
and similarly if $i>m$. However, for $i=m$ the similar deformation
$$
\bigl(t^m, (t-\epsilon)^m, t(t-\epsilon)\bigr),
$$
gives arcs that do not pass through the origin.
\medskip

These suggest that  $\sharcr(xy=z^m)$ should have $4m-6$
connected components.

\end{exmp}

The above example shows that the  method of Proposition \ref{sharcr.finite.prop}
can give too many candidates for connected components, when
applied to the minimal resolution. One can, however, hope that
the minimal dlt modification  is again the right choice.

\begin{prob}
Let $X$ be a normal, real algebraic surface with minimal dlt modification
$g^{\rm dlt}:\bigl(E^{\rm dlt}\subset X^{\rm dlt}\bigr)\to X$. 
Does composing with $g^{\rm dlt} $ induce a bijection between
the (connected or irreducible)
 components of $\sharcr\bigl(E^{\rm dlt}\subset X^{\rm dlt}\bigr)$
 and  $\sharcr(0\in X)$?
\end{prob}

The proof of Theorem \ref{dlt.nonquot.thm}
relies on the fundamental group of the complex link of a singularity;
we know no analogs of it in the real case.

%\bibliography{refs-main/refs}
\def\cprime{$'$} \def\cprime{$'$} \def\cprime{$'$} \def\cprime{$'$}
  \def\cprime{$'$} \def\cprime{$'$} \def\dbar{\leavevmode\hbox to
  0pt{\hskip.2ex \accent"16\hss}d} \def\cprime{$'$} \def\cprime{$'$}
  \def\polhk#1{\setbox0=\hbox{#1}{\ooalign{\hidewidth
  \lower1.5ex\hbox{`}\hidewidth\crcr\unhbox0}}} \def\cprime{$'$}
  \def\cprime{$'$} \def\cprime{$'$} \def\cprime{$'$}
  \def\polhk#1{\setbox0=\hbox{#1}{\ooalign{\hidewidth
  \lower1.5ex\hbox{`}\hidewidth\crcr\unhbox0}}} \def\cdprime{$''$}
  \def\cprime{$'$} \def\cprime{$'$} \def\cprime{$'$} \def\cprime{$'$}
\providecommand{\bysame}{\leavevmode\hbox to3em{\hrulefill}\thinspace}
\providecommand{\MR}{\relax\ifhmode\unskip\space\fi MR }
% \MRhref is called by the amsart/book/proc definition of \MR.
\providecommand{\MRhref}[2]{%
  \href{http://www.ams.org/mathscinet-getitem?mr=#1}{#2}
}
\providecommand{\href}[2]{#2}

\medskip

\noindent Princeton University, Princeton NJ 08544-1000

{\begin{verbatim}kollar@math.princeton.edu\end{verbatim}}
\medskip

\noindent R\'enyi Institute of Mathematics, Budapest, Hungary

{\begin{verbatim}nemethi.andras@renyi.mta.hu\end{verbatim}}

\end{document}